\documentclass[12pt]{amsart}
\usepackage{amssymb}
\usepackage{amsfonts}
\usepackage{latexsym}
\usepackage{amscd}
\usepackage[mathscr]{euscript}
\usepackage{xy} \xyoption{all}

\newcommand{\N}{{\mathbb{N}}}
\newcommand{\Z}{{\mathbb{Z}}}

\newcommand{\uloopr}[1]{\ar@'{@+{[0,0]+(-4,5)}@+{[0,0]+(0,10)}@+{[0,0] +(4,5)}}^{#1}}
\newcommand{\uloopd}[1]{\ar@'{@+{[0,0]+(5,4)}@+{[0,0]+(10,0)}@+{[0,0]+ (5,-4)}}^{#1}}
\newcommand{\dloopr}[1]{\ar@'{@+{[0,0]+(-4,-5)}@+{[0,0]+(0,-10)}@+{[0, 0]+(4,-5)}}_{#1}}
\newcommand{\dloopd}[1]{\ar@'{@+{[0,0]+(-5,4)}@+{[0,0]+(-10,0)}@+{[0,0 ]+(-5,-4)}}_{#1}}

\newcommand{\luloop}[1]{\ar@'{@+{[0,0]+(-8,2)}@+{[0,0]+(-10,10)}@+{[0, 0]+(2,2)}}^{#1}}

\newtheorem{lem}{Lemma}[section]
\newtheorem{corol}[lem]{Corollary}
\newtheorem{theor}[lem]{Theorem}
\newtheorem{prop}[lem]{Proposition}

\newtheorem{defi}[lem]{Definition}
\newtheorem{defis}[lem]{Definitions}
\newtheorem{exem}[lem]{Example}

\newtheorem{point}{}

\begin{document}
\title[Classification of Leavitt path algebras]{The classification question for Leavitt path algebras}
\author{G. Abrams}
\address{Department of Mathematics, University of Colorado,
Colorado Springs CO 80933 U.S.A.} \email{abrams@math.uccs.edu}
\urladdr{http://mathweb.uccs.edu/gabrams/}
\author{P. N. \'{A}nh}
\address{R\'enyi Institute of Mathematics, Hungarian Academy of
Sciences, 1364 Budapest, Pf. 127 Hungary} \email{anh@renyi.hu}
\author{A. Louly}\author{E. Pardo}
\address{Departamento de Matem\'aticas, Universidad de C\'adiz,
Apartado 40, 11510 Puerto Real (C\'adiz), Spain.}\email{louly.adel@uca.es, enrique.pardo@uca.es}
\urladdr{http://www2.uca.es/dept/matematicas/PPersonales/PardoEspino/index.HTML}
\thanks{The second author is supported partly by Hungarian National Foundation
for Scientific Research grant no. K61007. During Fall 2006 he was also supported partly by The
Colorado College, the University of Colorado at Colorado Springs, and Professor Michael Siddoway.
The third author was partially supported by a Post-doctoral fellow of the A.E.C.I. (Spain). The
fourth author was partially supported by the DGI and European Regional Development Fund, jointly,
through Projects MTM2004-00149 and MTM2007-60338, by PAI III grants FQM-298 and P06-FQM-1889 of the
Junta de Andaluc\'{\i}a, and by the Comissionat per Universitats i Recerca de la Generalitat de
Catalunya.} 
%
\dedicatory{Dedicated to the memory of Manuel Berrocoso Senior}
%
\begin{abstract}
We prove an algebraic version of the Gauge-Invariant Uniqueness Theorem, a result which gives
information about the injectivity of certain homomorphisms between ${\mathbb Z}$-graded algebras.
As our main application of this theorem, we obtain isomorphisms between the Leavitt path algebras
of specified graphs.  From these isomorphisms we are able to achieve two ends.  First, we show that
the $K_0$ groups of various sets of purely infinite simple Leavitt path algebras, together with the
position of the identity element in $K_0$, classify the algebras in these sets up to isomorphism.
Second, we show that the isomorphism between matrix rings over the classical Leavitt algebras,
established previously using number-theoretic methods,  can be reobtained via appropriate
isomorphisms between Leavitt path algebras.
\end{abstract}

\maketitle

\section*{Introduction}\label{Introduction}

Throughout this article $E$ will denote a finite directed graph, and $K$ will denote an arbitrary
field.  The {\it Leavitt path algebra of $E$ with coefficients in $K$}, denoted $L_K(E)$, has
received significant attention  over the past few years, both from algebraists as well as from
analysts working in operator theory. (The precise definition of $L_K(E)$ is given below.)   When
$K$ is the field ${\mathbb C}$ of complex numbers, the algebra $L_K(E)$ has exhibited surprising
similarity to its C$^*$-algebra counterpart C$^*(E)$, the {\it Cuntz-Krieger graph} C$^*$-{\it
algebra of} $E$. In part motivated by the Gauge-Invariant Uniqueness Theorem of operator theory, we
prove in Section 1 of this article a result which guarantees the injectivity of certain algebra
homomorphisms based on a specified action of the field as automorphisms on the codomain. With this
result in hand, we use it (and other results about homomorphisms from Leavitt path algebras
established here) in Section \ref{XXX} to produce isomorphisms between Leavitt path algebras.
Specifically, we show in Theorems \ref{meu1_gen} and \ref{outsplitiso} how, starting with the graph
$E$ and specified configurations of vertices and edges in $E$, to explicitly construct graphs $E'$
having $L_K(E)\cong L_K(E')$.   In Section \ref{purelyinfinitesimple} we apply these isomorphisms
to obtain Proposition \ref{onlyapath}, a result about isomorphisms between purely infinite simple
unital Leavitt path algebras. Subsequently, in Section \ref{graphs} we use these isomorphisms to
answer specific cases of The Classification Question for purely infinite simple unital Leavitt path
algebras.  We establish Propositions \ref{twovertices} and \ref{threevertices}, the algebraic
counterparts of specific pieces of the Kirchberg - Phillips Theorem of C$^*$-algebras (see
\cite[Section 3]{T1} for a description). These two results establish, respectively, that if $E$ and
$F$ are finite directed graphs having two (respectively three) vertices and no parallel edges, and
the Leavitt path algebras $L_K(E)$ and $L_K(F)$ are purely infinite simple, then $L_K(E)\cong
L_K(F)$ if and only if the Grothendieck groups $K_0(L_K(E))$ and $K_0(L_K(F))$ are isomorphic via
an isomorphism which takes $[1_{L_K(E)}]$ to $[1_{L_K(F)}]$. We close the article by showing how
\cite[Theorem 4.14]{AbAnhP} may be reestablished using the results and techniques of this article.

We briefly recall some graph-theoretic definitions and properties; more complete explanations and
descriptions can be found in \cite{AA1}. A \emph{graph} $E=(E^0,E^1,r,s)$ consists of two countable
sets $E^0,E^1$ and maps $r,s:E^1 \to E^0$.  (Some authors use the phrase `directed' graph for this
structure.)  The elements of $E^0$ are called \emph{vertices} and the elements of $E^1$
\emph{edges}. (We emphasize that loops and multiple / parallel edges are allowed.) If $s^{-1}(v)$
is a finite set for every $v\in E^0$, then the graph is called \emph{row-finite}. A vertex $v$ for
which $s^{-1}(v)$ is empty is called a \emph{sink}, while a vertex $w$ for which $r^{-1}(w)$ is
empty is called a \emph{source}. A \emph{path} $\mu$ in a graph $E$ is a sequence of edges
$\mu=e_1\dots e_n$ such that $r(e_i)=s(e_{i+1})$ for $i=1,\dots,n-1$. In this case,
$s(\mu):=s(e_1)$ is the \emph{source} of $\mu$, $r(\mu):=r(e_n)$ is the \emph{range} of $\mu$, and
$n$ is the \emph{length} of $\mu$. An edge $e$ is an {\it exit} for a path $\mu = e_1 \dots e_n$ if
there exists $i$ such that $s(e)=s(e_i)$ and $e\neq e_i$. If $\mu$ is a path in $E$, and if
$v=s(\mu)=r(\mu)$, then $\mu$ is called a \emph{closed path based at $v$}. If $\mu= e_1 \dots e_n$
is a closed path based at $v = s(\mu)$ and $s(e_i)\neq s(e_j)$ for every $i\neq j$, then $\mu$ is
called a \emph{cycle}. We say that a graph $E$ satisfies \emph{Condition} (L) if every cycle in $E$
has an exit. For $n\ge 2$ we define $E^n$ to be the set of paths of length $n$, and
$E^*=\bigcup_{n\ge 0} E^n$ the set of all paths.

 The following notation is standard.
  Let $A$ be a $p\times p$ matrix having non-negative integer entries (i.e., $A = (a_{ij})\in {\rm M}_p(\Z^+)$).  The
  graph $E_A$ is defined by setting $(E_A)^0=\{v_1,v_2,...,v_p\}$, and defining $(E_A)^1$ by inserting  exactly $a_{ij}$ edges in $E_A$ having source vertex $v_i$ and
  range vertex $v_j$.   Conversely, if $E$ is a finite graph with vertices $\{v_1,v_2,...,v_p\}$, then we define the
   {\it incidence matrix} $A_E$ {\it of} $E$ by setting
  $(A_E)_{ij}$ as the number of edges in $E$ having source vertex $v_i$ and range vertex $v_j$.


\begin{defi}\label{definition}  {\rm Let $E$ be any row-finite graph, and $K$ any field. The
{\em Leavitt path $K$-algebra} $L_K(E)$ {\em of $E$ with coefficients in $K$} is  the $K$-algebra
generated by a set $\{v\mid v\in E^0\}$ of pairwise orthogonal idempotents, together with a set of
variables $\{e,e^*\mid e\in E^1\}$, which satisfy the following relations:

(1) $s(e)e=er(e)=e$ for all $e\in E^1$.

(2) $r(e)e^*=e^*s(e)=e^*$ for all $e\in E^1$.

(3) $e^*e'=\delta _{e,e'}r(e)$ for all $e,e'\in E^1$.

(4) $v=\sum _{\{ e\in E^1\mid s(e)=v \}}ee^*$ for every vertex $v\in E^0$ for which $s^{-1}(v)$ is
nonempty.

}

\end{defi}

When the role of the coefficient field $K$ is not central to the discussion, we will often denote
$L_K(E)$ simply by $L(E)$.

The set $\{e^*\mid e\in E^1\}$ will be denoted by $(E^1)^*$. We let $r(e^*)$ denote $s(e)$, and we
let $s(e^*)$ denote $r(e)$. If $\mu = e_1 \dots e_n$ is a path, then we denote by $\mu^*$ the
element $e_n^* \dots e_1^*$ of $L_K(E)$.

For any subset $H$ of $E^0$, we will denote by $I(H)$ the ideal of $L_K(E)$ generated by $H$.

An alternate description of $L_K(E)$ is given in \cite{AA1}, where it is described in terms of a
free associative algebra modulo the appropriate relations indicated in Definition \ref{definition}
above. As a consequence, if $A$ is any $K$-algebra which contains elements satisfying these same
relations, then there is a (unique) $K$-algebra homomorphism from $L_K(E)$ to $A$ mapping the
generators of $L_K(E)$ to their appropriate counterparts in $A$.  We will refer to this conclusion
as the Universal Homomorphism Property of $L_K(E)$.   See also \cite[Remark 2.5]{T2}.

Many well-known algebras arise as the Leavitt path algebra of a row-finite graph. For instance, the
classical Leavitt algebras $L_n$ for $n\ge 2$ (see Definition \ref{Leavittalgebra} below) arise as
the algebras $L(R_n)$ where $R_n$ is the ``rose with $n$ petals" graph
$$
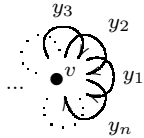
\xymatrix{ & {\bullet^v} \ar@(ur,dr) ^{y_1} \ar@(u,r) ^{y_2} \ar@(ul,ur) ^{y_3} \ar@{.} @(l,u) \ar@{.} @(dr,dl)
\ar@(r,d) ^{y_n} \ar@{}[l] ^{\ldots} }$$ The full $n\times n$ matrix algebra over $K$ arises as the
Leavitt path algebra of the oriented $n$-line graph
$$\xymatrix{{\bullet}^{v_1} \ar [r] ^{e_1} & {\bullet}^{v_2} \ar [r] ^{e_2} & {\bullet}^{v_3} \ar@{.}[r] &
{\bullet}^{v_{n-1}} \ar [r] ^{e_{n-1}} & {\bullet}^{v_n}} $$ while the Laurent polynomial algebra
$K[x,x^{-1}]$ arises as the Leavitt path algebra of the ``one vertex, one loop" graph
$$
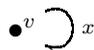
\xymatrix{{\bullet}^{v} \ar@(ur,dr) ^x}$$ Constructions such as direct sums and the formation of matrix rings produce
additional examples of Leavitt path algebras.

If $E$ is a finite graph then $L_K(E)$ is unital, with $\sum _{v\in E^0} v=1_{L_K(E)}$. Conversely,
if $L_K(E)$ is unital, then $E^0$ is finite.
 If $E^0$ is infinite then $L_K(E)$
is a ring with a set of local units; one such set of local units consists of sums of distinct
elements of $E^0$.  $L_K(E)=\bigoplus_{n\in \Z} L_K(E)_n$ is a ${\mathbb Z}$-graded $K$-algebra,
spanned as a $K$-vector space by $\{pq^* \mid p,q$ are paths in $E\}$.  In particular, for each
$n\in \mathbb{Z}$, the {\it degree} $n$ component $L_K(E)_n$ is spanned by elements of the form
$\{pq^* \mid {\rm length}(p)-{\rm length}(q)=n\}$.  The degree of an element $x$, denoted $deg(x)$,
is the minimum integer $n$ for which $x\in \bigoplus_{m\leq n} L_K(E)_m$. The set of
\emph{homogeneous elements} is $\bigcup_{n\in {\mathbb Z}} L_K(E)_n$, and an element of $L_K(E)_n$
is said to be $n$-\emph{homogeneous} or \emph{homogeneous of degree} $n$. The $K$-linear extension
of the assignment $pq^* \mapsto qp^*$ (for $p,q$ paths in $E$) yields an involution on $L_K(E)$,
which we denote simply as ${}^*$.

Information regarding the ``C$^*$-algebra of a graph", also known as the ``Cuntz-Krieger graph
C$^*$-algebra", may be found in \cite{R}.

We thank the referee for a very careful and helpful review of this article.

\section{Injectivity of algebra maps}\label{Injectivity}

Our central theme in this article is a description of isomorphisms between Leavitt path algebras.
As we shall see, oftentimes we encounter a situation in which we have defined a surjective ring
homomorphism between two such algebras, and seek to determine whether the map is injective.  The
main result of this section, Theorem \ref{AGIUT} (which we refer to as the Algebraic
Gauge-Invariant Uniqueness Theorem or AGIUT for short), provides a tool for doing just that. The
AGIUT is a consequence of results for general ${\mathbb Z}$-graded $K$-algebras.

In fact, there are many results aside from the AGIUT which provide similar tools by which we can
establish the injectivity of various algebra homomorphisms. We present two such results in the
following lemmas.

\begin{lem}\label{gradmap1}
Let $E$ be a row-finite graph, let $A$ be a $\Z$-graded $K$-algebra, and let $f:L_K(E)\rightarrow
A$ be an algebra map such that $f(v)\ne 0$ for every $v\in E^0$. If $f$ is graded, then $f$ is
injective.
\end{lem}
\begin{proof}
Since $f$ is a graded map, $\mbox{Ker}(f)$ is a graded ideal of $L_K(E)$. By \cite[Theorem
5.3]{AMFP}, there exists a subset $X$ of $E^0$ such that $\mbox{Ker}(f)=I(X)$. Since $f(v)\ne 0$
for every $v\in E^0$, we get $X=\emptyset$, whence $\mbox{Ker}(f)=0$ as desired.
\end{proof}

The method used in the previous proof is to guarantee that $\mbox{Ker}(f)$ does not contain any
vertices. The proof of the following lemma uses  a similar line of reasoning.

\begin{lem}\label{gradmap2}
Let $E$ be a row-finite graph satisfying Condition  (L), let $A$ be a $K$-algebra, and let
$f:L_K(E)\rightarrow A$ be an algebra map. If $f(v)\ne 0$ for every $v\in E^0$, then $f$ is
injective.
\end{lem}
\begin{proof}
By \cite[Lemma 3.9]{AA1}, Condition (L) yields that $J\cap E^0\neq\emptyset$ for every nonzero
ideal $J$ of $L_K(E)$. Since $f(v)\ne 0$ for every $v\in E^0$, we conclude that $\mbox{Ker}(f)=0$
as desired.
\end{proof}

With the hypotheses of the previous two lemmas in mind, we seek an injectivity result in situations
in which the map is not graded, and the graph contains cycles with no exits.  Such a result is the
essence of the AGIUT (Theorem \ref{AGIUT}).

\begin{defis}\label{general_standard_action}
{\rm Let $K$ be a field, and let $A$ be a ${\Z}$-graded algebra over $K$. For $t \in K^* = K
\setminus  \{0\}$  and   $a$   any homogeneous element of $A$ of degree $d$, set
$$\tau_t(a) = t^d a,$$
 and extend $\tau_t$ to
all of $A$ by linearity. It is easy to show that $\tau _t$ is a $K$-algebra automorphism of $A$ for
each $t\in K^*$. Then $\tau : K^* \rightarrow {\rm Aut}_K(A)$ is an action of $K$ on $A$, which we
call the {\it gauge action} of $K$ on $A$.

If $I$ is an ideal of $A$, we say that $I$ is {\it gauge-invariant} in case $\tau_t(I) =I$ for each
$t \in K^*$. This condition is equivalent to requiring that $\tau_t(I) \subseteq  I$ for every $t
\in K^*$, since $\tau_{t^{-1}}(I) \subseteq I$ gives $I \subseteq \tau_t(I)$.}
\end{defis}

The previous definition of the gauge action draws its motivation as follows.  Let $A$ be a
$\Z$-graded $K$-algebra (e.g., a Leavitt path algebra $L_K(E)$) which is generated by homogeneous
elements of degree 1 and $-1$. Then the multiplicative group $K^*$ acts naturally on $A$ by sending
elements $a$ of degree $\epsilon$, where $\epsilon= 1, -1$ to $t^{\epsilon}a$ for each nonzero
element $t\in K$. In particular if $A$ is an involutorial $K-$algebra over a field $K$ with
involution, then any homogeneous element of degree $-1$ is the image of such an element of degree
1, and hence the unitary group of $K$ (i.e., the group of $t\in K$ with $tt^*=1$) acts naturally on
$A$. This natural action of the unit circle on the Cuntz algebra ${\mathcal O}_n$ is in part the
motivation for our description of an algebraic gauge action of $K$ on $A$. (For additional
information see \cite[p. 198]{Pim}.)

The next result establishes a  relationship between graded and gauge-invariant ideals of any
$\Z$-graded algebra.

\begin{prop}\label{generalgradedrings} Let $K$ be a field, let $A$
be a ${\Z}$-graded $K$-algebra, and let $I$  be an ideal of  $A$. Let $\tau : K^* \rightarrow {\rm
Aut}_K(A)$ be the gauge action of $K$ on $A$.
\begin{enumerate}
\item If $I$ is generated as an ideal of $A$ by elements of degree
0, then $I$ is gauge-invariant.\item If $K$ is infinite, and if $I$ is gauge-invariant, then $I$ is
graded.
\end{enumerate}
\end{prop}
\begin{proof}
Statement (1) is clear, as $\tau_t$ fixes the degree zero elements of $A$ for each $t \in K^*$.

For statement (2), we prove the contrapositive. So suppose $I$ is not graded. We seek to show that
$I$ is not gauge-invariant.   For each $a\in A$ let $h(a)$ denote the number of nonzero homogeneous
graded components of $a$. Since $I$ is not graded there exists an element $a \in I$ for which, in
the decomposition $a = \sum a_j$ into a sum of its homogeneous components, at least two of the
$a_j$  are not in $I$. Let $T \subseteq I$ denote those elements of $I$ which, when written in
homogeneous decomposition,  have the property that no nonzero homogeneous component is in $I$.
Since $I$ is not graded, $T \neq \emptyset$. Let  $b\in T$ such that $h(b)=\min\{h(t)|t\in T\}$.
Note that $h(b) \geq 2$. Let $m,n \in \Z$ for which the homogeneous components $b_m$ and $b_n$ are
each nonzero; assume without loss of
 generality that $n<m$. Because $K$ is infinite, we can find  $t \in K^*$ such that $t^m \neq t^n$. (Otherwise, we would have $t^m=t^n$ for all $t\in K^*$, so that every
element of $K^*$ would be a zero of $x^{m-n}-1 \in K[x]$, but such cannot happen in an infinite
field.)

To show that $I$ is not gauge-invariant, it suffices to show $\tau_t(I)\nsubseteq I$.  By
contradiction, assume $\tau_t (I)\subseteq I$, so that in particular $\tau_t (b) \in I$. We observe
that $b\in I$ gives $t^m  b \in I$, so $\tau_t (b) - t^m  b \in I$; we denote $\tau_t (b) - t^m b$
by $c$. Note that for each $i\in {\mathbb Z}$, the $i$-component of $c$ is $c_i = (t^i - t^m)b_i$.
Thus we have $c_m=0$, but $c_n = (t^n - t^m)b_n \neq 0$ (and so in particular $c\neq 0$). But $c_m
= 0$ gives $h(c) < h(b)$, so, by minimality, at least one of the nonzero components of $c$ is in
$I$. That is, for some $p \in \Z$, $(t^p - t^m)b_p$ is a nonzero element of $I$. But then $b_p$ is
a nonzero element of $I$, which contradicts our choice of $b$.
\end{proof}

We now apply this result in the context of Leavitt path algebras. For clarity, we present here the
definition of the gauge action of $K$ on the Leavitt path algebra $L_K(E)$ of the row-finite graph
$E$.

\begin{defi}\label{standard_action}
{\rm Let $E$ be a row-finite graph, and let $K$ be a field. Then the {\it gauge action} $\tau$ {\it
of} $K$ {\it on}
  the Leavitt path algebra $L_K(E)$ (denoted
sometimes by $\tau^E$ for clarity) is given by
$$\begin{array}{cccc}
  \tau^E : & K^* & \rightarrow & \mbox{Aut}_K(L_K(E)) \\
   & t & \mapsto & \tau^E _t \\
\end{array}$$
as follows: for every $t\in K^*$, for every $v\in E^0$, and for every $e\in E^1$
$$\begin{array}{cccc}
  \tau^E _t: & L_K(E) & \rightarrow & L_K(E) \\
   & v & \mapsto & v \\
   & e & \mapsto & t e \\
   & e^* & \mapsto & t^{-1} e^* \\
\end{array}$$
and then extend linearly and multiplicatively to all of $L_K(E)$.}
\end{defi}

For a graph $E$, the set of graded ideals of $A = L_K(E)$ is denoted by $\mathcal{L}_{\text{gr}}$.

\begin{prop}\label{graded}
Let $E$ be a row-finite graph, let $K$ be an infinite field, and let $I$ be an ideal of $L_K(E)$.
Then $I\in \mathcal{L}_{\text{gr}}$ if and only if $I$ is gauge-invariant.
\end{prop}
\begin{proof}
If $I\in \mathcal{L}_{\text{gr}}$, then $I=I(H)$ for some $H \subseteq E^0$ by \cite[Theorem
5.3]{AMFP}. Thus $I$ is generated by elements of degree zero, and so Proposition
\ref{generalgradedrings}(1) applies. The converse follows immediately from Proposition
\ref{generalgradedrings}(2).
\end{proof}

We note that the implication $I\in \mathcal{L}_{\text{gr}}$ implies $I$ is gauge-invariant holds
for any field $K$, finite or infinite. In contrast, we now show that the converse implication of
Proposition \ref{graded} is never true for any finite field.

\begin{prop}\label{falseiffinite}
For any finite field $K$ there exists a graph $E$ such that the Leavitt path algebra $L_K(E)$
contains a non-graded ideal which is gauge-invariant.
\end{prop}
\begin{proof}
If we denote card($K$) by $m + 1$, then $t^m = 1$ for all $t \in K^*$. Let $E$ be the graph
$$\xymatrix{{\bullet}^{v} \ar@(ur,dr) ^x}$$
so that, as noted previously, $L_K(E)\cong K[x,x^{-1}]$.  In particular we have $\tau_t(1+ x^m) =
1+ x^m$ for all $t \in K^*$. This then yields that the ideal $I = <1+ x^m>$ of $L_K(E)$ is
gauge-invariant.  But it is well known (or it can be shown using an argument similar to that given
in the proof of \cite[Theorem 3.11]{AA1}) that $I$ is not a graded ideal of $K[x,x^{-1}]$.
\end{proof}

We are now in position to present the main application of these ideas.

\begin{theor}\label{AGIUT}
(The Algebraic Gauge-Invariant Uniqueness Theorem.)  Let $E$ be a row-finite graph, let $K$ be an
infinite field, and let $A$ be a $K$-algebra. Suppose
$$\phi: L_K(E)\rightarrow A$$
 is a $K$-algebra homomorphism such that $\phi (v)\ne 0$ for every $v\in E^0$. If there exists a group
action $\sigma : K^* \rightarrow \mbox{Aut}_K(A)$ such that $\phi \circ \tau^E _t=\sigma _t \circ
\phi$ for every $t\in K^*$,  then $\phi $ is injective.
\end{theor}
\begin{proof}
Let $I=\mbox{Ker}(\phi)$. Then for every $a\in I$ and for every $t\in K^*$, $\phi (\tau^E
_t(a))=\sigma _t (\phi (a)) = \sigma_t(0) =0$, whence $\tau^E _t(a)\in \mbox{Ker}(\phi)=I$. Thus
for every $t\in K^*$ we have $\tau^E _t(I)\subseteq I$, so that $I$ is gauge-invariant.   Hence
$I\in \mathcal{L}_{\text{gr}}$ by Proposition \ref{graded}. In particular, if $I \ne \{0\}$, then
$I\cap E^0\ne \emptyset$ by \cite[Proposition 5.2 and Theorem 5.3]{AMFP}, contradicting the
hypothesis that $\phi (v)\ne 0$ for every $v\in E^0$.
\end{proof}

In both \cite{AA3} and \cite{T2} an analysis of Leavitt path algebras for non-row-finite graphs is
carried out.  We conclude Section \ref{Injectivity}  by noting that all the results (and their
proofs) presented in this section hold verbatim in this more general not-necessarily-row-finite
setting.    In particular, Lemma \ref{gradmap1} generalizes as \cite[Theorem 4.8]{T2}, while Lemma
\ref{gradmap2} generalizes as \cite[Theorem 6.8]{T2}.

\section{Isomorphisms: general results}\label{XXX}

In this section we will apply the results of Section 1 to draw conclusions about isomorphisms
between Leavitt path algebras.
 The main goal in establishing such isomorphisms is as follows:  starting with a
graph $E$, we seek a systematic method to produce various graphs $F$ for which $L_K(E)\cong
L_K(F)$. As such, we refer to our
 two main results (Theorems \ref{meu1_gen} and \ref{outsplitiso}) as ``Change the Graph" theorems.    These results
  in turn will allow us to verify that a specific set of
Leavitt path algebras is determined up to isomorphism by $K_0$ data.

In our first such result, we show how to ``bundle" specific sets of edges, and subsequently replace
the bundled sets by a single edge.

\begin{defi}\label{transfer_rest}
{\rm Let $E$ be a row-finite graph, and let $v\ne w\in E^0$ be vertices which are not sinks. If
there exists an
 injective map
$$\theta :s^{-1}(w)\rightarrow s^{-1}(v)$$
such that $r(e)=r(\theta (e))$ for every $e\in s^{-1}(w)$, we define the {\it shift graph} from $v$
to $w$, denoted
$$F=E(w\hookrightarrow v),$$
 as follows:
\begin{enumerate}
\item $F^0=E^0$.\item $F^1=(E^1\setminus \theta(s^{-1}(w)))\cup \{
f_{v,w}\}$, where $f_{v,w}\not\in E^1$, $s(f_{v,w})=v$ and $r(f_{v,w})=w$.
\end{enumerate}
Although the definition the graph $F=E(w\hookrightarrow v)$ depends on the map $\theta$,
 in order to make the notation less cumbersome we suppress $\theta$ in the notation.  This will
 cause no confusion throughout the sequel.}
\end{defi}

\begin{exem}\label{shiftexample}
{\rm Consider the following graphs:
$$
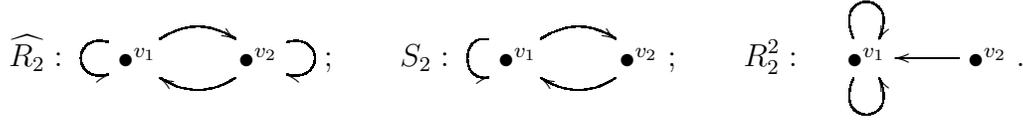
\widehat{R_2}: {
\def\labelstyle{\displaystyle}
\xymatrix{\bullet^{v_1}\dloopd{} \ar@/^1pc/ [r] & {\bullet}^{v_2} \ar@/^1pc/ [l] \uloopd{} }};
\qquad {
\def\labelstyle{\displaystyle}S_2:\quad \xymatrix{
{\bullet}^{v_1} \ar@(ul,dl)  \ar@/^1pc/ [r] & {\bullet}^{v_2} \ar@/^1pc/ [l]  }}; \qquad
{R_2^2}:\quad {
\def\labelstyle{\displaystyle}
\xymatrix{ \bullet^{v_1}\uloopr{}\dloopr{} & \bullet^{v_2} \ar[l]}}.
$$

Then notice that $S_2= \widehat{R_2}(v_1\hookrightarrow v_2)$ and $S_2=R_2^2(v_2\hookrightarrow
v_1)$. }
\end{exem}

Recall that a graph $E$ satisfies Condition (L) in case every cycle in $E$ has an exit.   It is
clear that if $E$ satisfies Condition (L), then so also does $E(w\hookrightarrow v)$ for any shift
graph constructed from $E$.

We are now in position to prove the first of two ``Change the Graph" Theorems.

\begin{theor}\label{meu1_gen}
Let $E$ be a row-finite graph, and let $v\ne w\in E^0$ be vertices which are not sinks. If there
exists an injection
$$\theta :s^{-1}(w)\rightarrow s^{-1}(v)$$
such that $r(e)=r(\theta (e))$ for every $e\in s^{-1}(w)$, then $L(E(w\hookrightarrow v))$ is a
homomorphic image of $L(E)$. Moreover, if either:
\begin{enumerate}
\item $E$ satisfies Condition {\rm (L)}, or
\item the field $K$ is infinite,
\end{enumerate}
then there exists a $K$-algebra isomorphism $\varphi: L(E) \rightarrow L(E(w\hookrightarrow v))$.
(The isomorphism $\varphi$
 is
 not an isomorphism of $\Z$-graded $K$-algebras.)
\end{theor}
\begin{proof}
Let $F=E(w\hookrightarrow v)$, and let $s_E^{-1}(w)=\{ e_1, \dots ,e_n\}$. Given any $e_i\in
s_E^{-1}(w)$, we define in $L(F)$ the element
$$T_{e_i}=f_{v,w}e_i.$$
 Notice that $T_{e_i}\ne
0$ for every $e_i\in s_E^{-1}(w)$, and that $T_{e_i}\ne T_{e_j}$ whenever $i\ne j$. Now consider
the subalgebra $A$ of $L(F)$ generated by
$$\{ v, e, e^*, T_{e_i}, T_{e_i}^*\mid v\in E^0, e\in E^1\setminus \theta(s_E^{-1}(w)), {e_i}\in s_E^{-1}(w)\}.$$
Then, if $i\ne j$, we have
$$T_{e_i}^*T_{e_j}=(f_{v,w}e_i)^*(f_{v,w}e_j)=e_i^*f_{v,w}^*f_{v,w}e_j=
e_i^*e_j=0=\theta(e_i)^*\theta(e_j),$$ while $T_{e_i}^*T_{e_i}=e_i^*e_i=r(e_i)=r(\theta
(e_i))=\theta (e_i)^*\theta (e_i)$. Also, $s(T_{e_i})=s(f_{v,w})=v=s(\theta(e_i))$ and
$r(T_{e_i})=r(e_i)=r(\theta (e_i))$. Moreover, the only generators in $A$ starting in $v$ which do
not belong to $s^{-1}(v)\setminus \theta(s^{-1}(w))$ are of the form $T_{e_i}$ with $e_i\in
s_E^{-1}(w)$. Thus,
$$\sum\limits _{\{e\in s^{-1}(v)\setminus
\theta(s^{-1}(w))\}}{ee^*} +  \sum\limits _{i=1}^{n}T_{e_i}T_{e_i}^*= \sum\limits _{\{e\in
s^{-1}(v)\setminus \theta(s^{-1}(w))\}}{ee^*} + \sum\limits _{i=1}^{n}f_{v,w}e_i e_i^*f_{v,w}^*=$$
$$\sum\limits _{\{e\in s^{-1}(v)\setminus \theta(s^{-1}(w))\}}{ee^*} + f_{v,w}\left(\sum\limits
_{i=1}^{n}e_i e_i^*\right)f_{v,w}^*=\sum\limits _{\{e\in s^{-1}(v)\setminus
\theta(s^{-1}(w))\}}{ee^*} + f_{v,w}f_{v,w}^*=v.$$ Hence the generators of $A$ satisfy the same
relations as do the elements of the set $\{ v,e,e^*\mid v\in E^0, e\in E^1\}$ in $L(E)$. Thus by
the Universal Homomorphism Property of $L(E)$ there exists a unique algebra morphism extending the
natural bijection
$$\begin{array}{cccc}
  \varphi : & L(E) & \rightarrow & L(F) \\
   & e & \mapsto & e \\
   & g & \mapsto & T_{\theta^{-1}(g)} \\
   & e^* & \mapsto & e^* \\
   & g^* & \mapsto & T_{\theta^{-1}(g)}^* \\
   & v & \mapsto & v \\
\end{array}$$
for every $e\in E^1\setminus \theta(s_E^{-1}(w))$, every $g\in \theta(s_E^{-1}(w))$, and every
$v\in E^0$.

Since $e_i\in E^1\setminus \theta (s_E^{-1}(w))$ for every $e_i\in s_E^{-1}(w)$, we have $e_i\in
A$, whence
$$f_{v,w}=f_{v,w}w=f_{v,w}\sum\limits_{i=1}^{n}e_ie_i^*=
\sum\limits_{i=1}^{n}T_{e_i}e_i^*.$$ But $e_i\in s_E^{-1}(w)$ implies $T_{e_i}=T_{\theta
^{-1}(\theta(e_i))}=\varphi (\theta(e_i))$, so that $f_{v,w}= \sum\limits_{i=1}^{n}T_{e_i}e_i^*=
\varphi (\sum\limits_{i=1}^{n}\theta({e_i})e_i^*)$, and hence $\varphi$ is onto.

We note here that $\varphi$ is not a graded homomorphism, since $deg(g)=1$, while
$deg(\varphi(g))=deg(T_{\theta^{-1}(g)})=2.$ Thus Lemma \ref{gradmap1} does not apply in this
situation.

In the first case, if $E$ satisfies Condition (L), then the injectivity of $\varphi$ may be
established by Lemma \ref{gradmap2}.

For the second case, if $K$ is infinite, then for every $t\in K^*$ we can define the automorphism
$\alpha_t$ of $L_K(F)$ by the extension of
$$\begin{array}{cccc}
  \alpha _t : & L(F) & \rightarrow & L(F) \\
   & e & \mapsto & te \\
   & f_{v,w} & \mapsto & f_{v,w} \\
   & e^* & \mapsto & t^{-1}e^* \\
   & f_{v,w}^* & \mapsto & f_{v,w}^* \\
   & v & \mapsto & v \\
\end{array}$$
for every $e\in E^1\setminus \theta(s_E^{-1}(w))$ and every $v\in E^0$. In this way we get an
action $\alpha: K^*\rightarrow \mbox{Aut}_K(L_K(F))$.  It is straightforward to check that, for
every $t\in K^*$, $\varphi\circ \tau^F_t=\alpha _t\circ\varphi$, where $\tau_t^F$ is the gauge
action of $K^*$ on $L_K(F)$. Thus, the injectivity of $\varphi$ derives from the AGIUT (Theorem
\ref{AGIUT}).
\end{proof}

\begin{exem}\label{isoofshiftexample}
{\rm Recall the graphs in Example \ref{shiftexample}. On one side, $S_2=
\widehat{R_2}(v_1\hookrightarrow v_2)$, whence $L(\widehat{R_2})\cong L(S_2)$ by Theorem
\ref{meu1_gen}. On the other side, $S_2=R_2^2(v_2\hookrightarrow v_1)$, so that $L(S_2)\cong
L(R_2^2)$ again by Theorem \ref{meu1_gen}.  }
\end{exem}

Theorem \ref{meu1_gen} admits a corresponding statement in the context of Cuntz-Krieger graph
C*-algebras. As far as we know, no such analogous result has been established elsewhere in the
C$^*$-algebra literature.  We do so here.

\begin{corol}\label{meu1_gen-C*}
Let $E$ be a row-finite graph, and let $v\ne w\in E^0$ be vertices which are not sinks. If there
exists an injection
$$\theta :s^{-1}(w)\rightarrow s^{-1}(v)$$
such that $r(e)=r(\theta (e))$ for every $e\in s^{-1}(w)$, then $C^*(E)\cong C^*(E(w\hookrightarrow
v))$.
\end{corol}
\begin{proof}
We will follow the C*-algebra notation (see e.g. \cite{Flow}). Let $F=E(w\hookrightarrow v)$. Given
any $e\in s_E^{-1}(w)$, we define in $C^*(F)$ the element
$$T_e=s_{f_{v,w}}s_ee.$$
Notice that $T_e\ne 0$ for every $e\in s_E^{-1}(w)$, and that $T_e\ne T_f$ whenever $e\ne f\in
s_E^{-1}(w)$. Now consider the $C^*$-subalgebra $A$ of $C^*(F)$ generated by
$$\mathcal{S}=\{ p_v, s_e, T_g\mid v\in E^0, e\in E^1\setminus \theta(s_E^{-1}(w)), g\in s_E^{-1}(w)\}.$$

To simplify notation, let $s_E^{-1}(w)=\{ e_1, \dots ,e_n\}$. Then, the same argument as in the
proof of Theorem \ref{meu1_gen} shows that $\mathcal{S}$ is a Cuntz-Krieger $E$-family, whence
there exists a unique $C^*$-algebra morphism extending the natural bijection
$$\begin{array}{cccc}
  \varphi : & C^*(E) & \rightarrow & C^*(F) \\
   & s_e & \mapsto & s_e \\
   & s_g & \mapsto & T_{\theta^{-1}(g)} \\
   & p_v & \mapsto & p_v \\
\end{array}$$
for every $e\in E^1\setminus \theta(s_E^{-1}(w))$, every $g\in \theta(s_E^{-1}(w))$, and every
$v\in E^0$. The same argument as above shows that $\varphi$ is an onto map, while injectivity is a
consequence of Gauge-Invariant Uniqueness Theorem for graph C$^*$-algebras (see e.g. \cite[Theorem
2.2]{R}), applied to the $\mathbb{T}$-action on $C^*(F)$ defined by
$$\begin{array}{cccc}
  \alpha _z : & C^*(F) & \rightarrow & C^*(F) \\
   & s_e & \mapsto & zs_e \\
   & s_{f_{v,w}} & \mapsto & s_{f_{v,w}} \\
   & p_v & \mapsto & p_v \\
\end{array}$$
for every $z\in \mathbb{T}$.
\end{proof}

In our second main result of this section, we show how to ``unbundle" specific sets of edges, and
subsequently replace these unbundled sets by a collection of new edges and new vertices. The
following definition is borrowed from \cite[Section 3]{Flow}.

\begin{defi}\label{outsplitdefi}
{\rm Let $E = ( E^0 , E^1 , r , s )$ be a row-finite graph. For each $v \in E^0$ which is not a
sink, partition $s^{-1} (v)$ into disjoint nonempty subsets $\mathcal{E}^1_v , \ldots ,
\mathcal{E}^{m(v)}_v$, where $m(v) \ge 1$.   (If $v$ is a sink, then we put $m(v)=0$.) Let
$\mathcal{P}$ denote the resulting partition of $E^1$. We form the {\em out-split graph $E_s (
\mathcal{P} )$ from $E$ using $\mathcal{P}$} as follows: Let
\begin{align*}
E_s ( \mathcal {P} )^0 &= \{ v^i \mid v \in E^0 , 1 \le i \le m(v) \} \cup \{ v :
m(v)=0 \} , \\
E_s ( \mathcal{P} )^1 &= \{ e^j \mid e \in E^1, 1 \le j \le m ( r (e) ) \} \cup \{ e : m ( r(e) ) =
0 \} ,
\end{align*}

\noindent and define $r_{E_s ( \mathcal{P} )} , s_{E_s ( \mathcal{P} )} : E_s ( \mathcal{P} )^1
\rightarrow E_s ( \mathcal{P} )^0$ for $e \in \mathcal{E}^i_{s(e)}$ by
\begin{align*}
s_{E_s ( \mathcal{P} )} ( e^j ) &= s(e)^i \text{ and } s_{E_s ( \mathcal{P} )} (e)
= s(e)^i  \\
r_{E_s ( \mathcal{P} )} ( e^j ) &= r(e)^j \text{ and }  r_{E_s ( \mathcal{P} )} ( e ) = r(e) .
\end{align*}}
\end{defi}

\begin{exem}\label{outsplitexample}
{\rm Consider the graph
$$
{R_2^2}:\quad {
\def\labelstyle{\displaystyle}
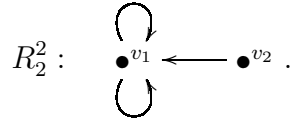
\xymatrix{ \bullet^{v_1}\uloopr{}\dloopr{} & \bullet^{v_2} \ar[l]}}.
$$
Let $\mathcal{P}$ be the partition of the edges of $R_2^2$ containing only one edge per subset.
Then the out-split graph of $R_2^2$ using $\mathcal{P}$ is}
$${
\def\labelstyle{\displaystyle}E:\quad 
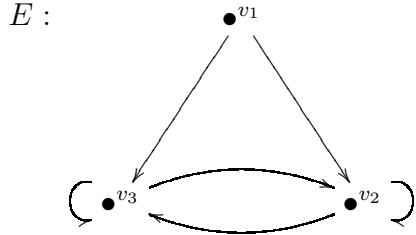
\xymatrix{ & \bullet^{v_1}
\ar[ddl] \ar[ddr]& \\ & & \\
{\bullet}^{v_3} \ar@(ul,dl) \ar@/^1pc/ [rr] & & {\bullet}^{v_2} \ar@(ur,dr) \ar@/^1pc/ [ll]  }}$$
\vspace{.2truecm}
\end{exem}

Similar to the graph C$^*$-algebra case, we get an isomorphism result for the Leavitt path algebras
of out-split graphs.  This result is the second
 of our two ``Change the Graph" theorems.

\begin{theor} \label{outsplitiso} {\rm (\cite[Theorem 3.2]{Flow})}
Let $E$ be a row-finite graph, $\mathcal{P}$ a partition of $E^1$ and $E_s ( \mathcal{P} )$ the
out-split graph formed from $E$ using $\mathcal{P}$. Then there is an
 isomorphism of $\Z$-graded $K$-algebras $\pi: L(E) \rightarrow L( E_s ( \mathcal{P} ) )$.
\end{theor}
\begin{proof}
The proof is essentially the same as that given in \cite[Theorem 3.2]{Flow}, except when showing
 the injectivity of the homomorphism. We include the argument here
for the sake of completeness.

Given $v \in E^0$ and $e \in E^1$, set $Q_v = v$ if $m(v)=0$, $T_e = e$ if $m(r(e))=0$,
\[
Q_v = \sum_{1 \le i \le m(v)} v^i \text{ if } m(v) \neq 0  \  \text{ and }  \  T_e = \sum_{1 \le j
\le  m ( r (e ) )} e^j \text{ if } m(r(e)) \neq 0 .
\]

\noindent Since $E$ is row-finite, all of these sums are finite. We claim that $\{ T_e , Q_v \mid e
\in E^1 , v \in E^0 \}$ is a family in $L( E_s ( \mathcal{P} ) )$ satisfying the same relations as
$\{v, e\mid v\in E^0, e\in E^1\}$.

The collection $\{Q_v | v\in E^0\}$ is a set of non-zero mutually orthogonal idempotents (since the
$Q_v$ are sums of idempotents satisfying the same properties). The elements $T_e$ for $e \in E^1$
clearly satisfy $T_e^*T_f=0$ whenever $e\ne f$, because they consist of sums of elements with the
same property. For $e \in E^1$ it is easy to see that $T_e^* T_e = Q_{r(e)}$.

For $e \in E^1$ with $m(r(e)) \neq 0$,  since $r_{E_s} ( e^j ) \neq r_{E_s} ( e^k )$, for $j \neq
k$, we have
\begin{equation} \label{reltwo}
T_e T_e^* = \Big( \sum_{1 \le j \le m ( r(e) )} e^j \Big) \Big( \sum_{1 \le k \le m ( r(e) )} e^k
\Big)^* = \sum_{1 \le j \le m(r(e))}{e^j}{e^j}^* .
\end{equation}

\noindent If $m(r(e))=0$ then $T_e T_e^* = e e^*$. For $v \in E^0$ and $1 \le i \le m(v)$ put
$$\mathcal{E}_{1,v}^i = \{ e \in
\mathcal{E}_v^i \mid m(r(e)) \geq 1 \} \mbox{ and } \mathcal{E}_{0,v}^i = \{ e \in \mathcal{E}_v^i
\mid m(r(e)) = 0 \}.$$
 If $v \in E^0$ is not a sink then $s^{-1}(v) = \bigcup_{i=1}^{m(v)}{\mathcal E}^i_v$ and for $1 \le i \le m(v)$
we have
\[
s_{E_s ( \mathcal{P} )}^{-1}
( v^i ) = \{ e^j \mid e \in \mathcal{E}_{1,v}^i , 1 \le j \le m(r(e)) \} \cup 
{\mathcal E}^i_{0,v}
\]

\noindent Hence using (\ref{reltwo}) we may compute
\begin{align*}
Q_v = \sum_{1 \le i \le m(v)} v^i &= \sum_{1 \le i \le m(v)} \sum_{e \in \mathcal{E}_{1,v}^i}
\sum_{1 \le j \le m(r(e))} e^j {e^j}^* + \sum_{1 \le i \le m(v)} \sum_{e \in
\mathcal{E}_{0,v}^i} e e^*  \\
&= \sum_{1 \le i \le m(v)} \; \sum_{e \in \mathcal{E}_v^i} T_e T_e^* = \sum_{\{e\mid s(e)=v\}} T_e
T_e^* ,
\end{align*}

\noindent completing the proof of our claim, since vertices $v \in E^0$ with $m(v)=0$ are sinks.

Then, by the Universal Homomorphism Property of $L(E)$ there is a homomorphism $\pi : L(E)
\rightarrow L ( E_s ( \mathcal{P} ))$ taking $e$ to $T_e$, $e^*$ to $T_e^*$ and $v$ to $Q_v$. To
prove that $\pi$ is onto we show that the generators of $L( E_s ( \mathcal{P} ) )$ lie in $L( T_e ,
Q_v )$, the subalgebra of $L(E_s(P))$ generated by $\{ T_e,T_e^*,Q_v\}$. Suppose that $w = v^j \in
E_s(\mathcal{P})^0$ is not a sink, set $e\in \mathcal{E}_v^j$, and pick $1\leq k\leq m(r(e))$. Then
$\{ f\in E_s(\mathcal{P})^1\mid s_{E_s(\mathcal{P})}(f)=v^j\}=\bigcup\limits _{e\in
\mathcal{E}_v^j}\{ e^k\mid 1\leq k\leq m(r(e))\}$, and we have
\[
v^j=\Big( \sum_{\{ f\in E_s(\mathcal{P})^1\mid s_{E_s(\mathcal{P})}(f)=v^j\}} f f^*  \Big)=\Big (
\sum\limits _{\{e\in \mathcal{E}_v^j\}}\sum\limits _{\{1\leq k\leq m(r(e))\}}e^k{e^k}^* \Big)=\Big
( \sum\limits _{\{e\in \mathcal{E}_v^j\}} T_eT_e^* \Big).
\]

If $w$ is a sink, then $w=Q_w$. Thus, $w \in L( T_e , Q_v )$.

If $e^j \in E_s ( \mathcal{P} )^1$ then $m(r(e)) \neq 0$. Since $r(e)^j \in L( T_e , Q_v )$ we have
$e^j = T_e r(e)^j
 \in L( T_e , Q_v )$.

If $e \in E_s ( \mathcal{P} )^1$ then $m(r(e)) = 0$ and so $e = T_e \in L( T_e , Q_v )$.

Since $Q_v$ is a sum of vertices and $T_e$ is a sum of edges, we get that $\pi$ is a $\Z$-graded
map, whence the injectivity of $\pi$ is guaranteed by Lemma \ref{gradmap1}, and the result follows.
\end{proof}

There is one specific partition which will play an important role throughout the sequel.

\begin{defi}\label{maximaloutsplitting}
{\rm For any row-finite graph $E$, the {\em maximal out-splitting} $\tilde{E}$ of $E$ is formed by
using the partition having $m(v)=|s^{-1}(v)|$ for every $v\in E^0$ which is not a sink.   In other
words, $\tilde{E}$ is the graph formed from $E$ by using the partition $\mathcal{P}$ of $E^1$ which
admits no refinements.}
\end{defi}

\begin{corol}\label{isototilde}
Let $E$ be a row-finite graph, and let $\tilde{E}$ denote the maximal out-splitting of $E$.  Then
$L_K(E)\cong L_K(\tilde{E})$ as
 $\Z$-graded $K$-algebras.
\end{corol}

As it turns out, the maximal out-splitting $\tilde{E}$ for a graph $E = ( E^0 , E^1 , r , s)$
without sinks is isomorphic to a graph which is well-known among graph theorists.  Recall that the
{\it dual graph} of a graph $E$ is the graph $\widehat{E} = ( E^1 , E^2 , r' , s' )$, where
$r'(ef)=f$ and $s'(ef)=e$.

\begin{prop}\label{maxoutsplitisotodual}
For any row-finite graph $E$ without sinks, the maximal out-splitting graph $\tilde{E}$ is
isomorphic to the dual graph $\widehat{E}$.
\end{prop}
\begin{proof}
Since the out-splitting is maximal, and $E$ is assumed to contain no sinks, we have
\[
\tilde{E}^0 = \{ v^e : s(e) = v \}  \ \text{ and } \ \tilde{E}^1 = \{ e^f : s(f)=r(e) \} .
\]

\noindent The maps $v^e \mapsto e$ and $e^f \mapsto ef$ are easily shown to induce an isomorphism
from $\tilde{E}$ to $\widehat{E}$.
\end{proof}

As a consequence of this proposition, it is reasonable to define the dual graph $\widehat{E}$ of
{\it any} row-finite graph $E$ to be its maximal out-splitting graph $\tilde{E}$. Thus, by
Corollary \ref{isototilde}, we get the following algebraic analog to a well-known result for graph
C*-algebras.

\begin{corol}\label{isomaxoutsplit}
If $E$ is any row-finite graph, then $L(E)\cong L(\widehat{E})$ as
 $\Z$-graded $K$-algebras.
\end{corol}

\section{The purely infinite simple unital case}\label{purelyinfinitesimple}

In this section we apply results from Section 2 to obtain information about the collection of
purely infinite simple unital Leavitt path algebras. Our first goal is to establish Corollary
\ref{simplyclass1}, which shows that, up to isomorphism, all purely infinite simple unital Leavitt
path algebras arise from a well-behaved subset of finite graphs. We start by reminding the reader
of the germane ring- and graph-theoretic ideas.

\begin{defis}\label{purelyinfinite}
{\rm Let $R$ be a ring.  A nonzero idempotent $e\in R$ is {\it infinite} if the right ideal $eR$
contains a proper direct summand isomorphic to $eR$. A ring $R$ is {\it purely infinite simple} if
\begin{enumerate}
\item $R$ is simple (i.e., $R$ contains no proper two-sided ideals), and\item every right ideal of $R$ contains an infinite
idempotent.
\end{enumerate}}
\end{defis}

\begin{defis}\label{cofinaldef}
 {\rm Let $E$ be a row-finite graph. If $v,w\in E^0$,  we say $w$  {\it connects to} $v$ if there is a path $\mu$ in
 $E$ for which $s(\mu)=w$ and $r(\mu)=v$.  If $c$ is a cycle in $E$, we say $w$ {\it connects to} $c$ if $w$ connects
 to some vertex $v$ in $c$.  A subset $H\subseteq E^0$ is
{\it hereditary} if whenever $w\in H$ and $v\in E^0$ and $w$ connects to $v$, then $v\in H$.  The
set $H$ is {\it saturated}
  if whenever $s^{-1}(v)\neq \emptyset$ and
$\{r(e):s(e)=v\}\subseteq H$, then $v\in H$.  The graph $E$ is called {\it cofinal} if the only
hereditary saturated subsets of $E^0$ are $\emptyset$ and $E^0$.   }
\end{defis}

The purely infinite simple Leavitt path algebras have been described in \cite[Theorem 11]{AA2}.
Specifically, $L(E)$ is purely infinite simple if and only if (i) $E$ is cofinal, (ii) $E$
satisfies Condition (L), and (iii) every vertex in $E^0$ connects to a cycle.

\begin{defi}
{\rm Given a row-finite graph $E$, we say that $E$ satisfies {\it Condition}  (Sing) if $E$
contains no parallel edges. Rephrased, $E$ satisfies Condition (Sing) if for every pair of vertices
$v,w\in E^0$, $\mbox{card}(\{ e\in E^1\mid s(e)=v \mbox{ and } r(e)=w\})\leq 1$.}
\end{defi}

\begin{prop}\label{onlyapath}
Let $E$ be a finite graph such that $L(E)$ is a purely infinite simple ring. If
$k=\mbox{card}(E^1)$, then for every $n\geq k$ there exists a graph $E_n$ such that:
\begin{enumerate}
\item $\mbox{card}(E_n^0)=n$, \item $E_n$ satisfies Condition
{\rm (Sing)}, and \item $L(E_n)\cong L(E)$.
\end{enumerate}
\end{prop}
\begin{proof}
By induction on $n$.   We start by establishing the result in the case $n=k$. By Corollary
\ref{isototilde}  we obtain $E_n=\tilde{E}$ having $\mbox{card}(E_n^0)=\mbox{card}(E^1)=k$ and
$L(E_n)\cong L(E)$. Since $E_n$ is the maximal out-splitting of $E$, it clearly satisfies Condition
(Sing).

Suppose that the result holds for some $n\geq k$; we will prove that it holds for $n+1$. So we pick
$E_n$ satisfying Condition  (Sing), $\mbox{card}(E_n^0)=n$ and $L(E_n)\cong L(E)$. Since $L(E_n)$
is purely infinite simple, $E_n$ contains at least one cycle having an exit by \cite[Theorem
11]{AA2}. Thus, there exists a vertex $v\in E_n^0$ such that $\mbox{card}(s^{-1}(v))\geq 2$.
Consider any nontrivial partition $\mathcal{P}$ with $s^{-1}(v)=\mathcal{E}_v^1\cup\mathcal{E}_v^2$
arbitrary, and for any $w\in E_n^0\setminus\{ v\}$ which is not a sink, let
$s^{-1}(w)=\mathcal{E}_w^1$ be the trivial partition. Then the out-split graph of $E_n$ by the
partition $\mathcal{P}$ satisfies $E_s(\mathcal{P})^0=(E_n^0\setminus \{ v\})\cup \{ v^1, v^2\}$,
so that $\mbox{card}(E_s(\mathcal{P})^0)=n+1$. Moreover, by Theorem \ref{outsplitiso},
$L(E_s(\mathcal{P}))\cong L(E_n)$. Also, as $E_s(\mathcal{P})$ is obtained by a partition of $E_n$,
it necessarily satisfies Condition (Sing). Hence by defining $E_{n+1}=E_s(\mathcal{P})$, the
induction step is established.
\end{proof}

Thus, in order to decide whether two purely infinite simple unital Leavitt path algebras $L(E)$ and
$L(F)$  are isomorphic, it is enough to consider the problem for isomorphic algebras $L(E_n)$ and
$L(F_n)$ where $|E_n^0| = |F_n^0|$ and each of $E_n,F_n$ satisfy Condition  (Sing). More formally,

\begin{corol}\label{simplyclass1}
An invariant $\mathcal{K}$ classifies purely infinite simple unital Leavitt path algebras up to
isomorphism if and only if, for each $n\in \N$,  $\mathcal{K}$ classifies up to isomorphism purely
infinite simple unital Leavitt path algebras of  graphs having $n$ vertices which satisfy Condition
{\rm (Sing)}.
\end{corol}

In fact, we can extend Proposition \ref{onlyapath} to stipulate that the new graphs have no
sources.

\begin{prop}\label{L:nosource}
Let $E$ be a finite graph such that $L(E)$ is a purely infinite simple ring. Then for any $n\geq
\vert E^1\vert$ there exists a finite graph $F$ such that:
\begin{enumerate}
\item $\vert F^0\vert = n$
\item $F$ satisfies Condition {\rm (Sing)}
\item $F$ has no sources, and
\item $L(E)\cong L(F)$.
\end{enumerate}
(We note that the isomorphism in (4) is not necessarily $\Z$-graded.)
\end{prop}
\begin{proof}
By Proposition \ref{onlyapath}, for any $n\geq \vert E^1\vert$, there exists a graph $E_n$ with $n$
vertices satisfying Condition (Sing) such that $L(E)\cong L(E_n)$. We show now that we can modify
$E_n$ if necessary to produce a graph $F$ for which $F$ has the desired properties.   Let $C$
denote the set of vertices of $E_n^0$ which lie in a closed simple path, and set $T=E_n^0\setminus
C$. Notice that if $T\ne \emptyset$, $v\in T$ and $v$ is not a source, then $s(r^{-1}(v))\subset
T$. For, suppose $w\in s(r^{-1}(v))\cap C$. Then, there exists a cycle $\mu$ such that $w\in \mu
^0$. By cofinality of $E_n$, there exists $\alpha \in E_n^*$ with $s(\alpha )=v$ and $r(\alpha )\in
\mu^0$. But then, $v\in C$, contradicting the assumption. Hence, if $S$ denotes the sources of
$E_n^0$, it is clear that every vertex in $T$ lies in the tree of some $x\in S$, and that the tree
of $S$ feeds into $C$. Since $E_n$ is finite, we can partition $T$ in layers as follows: $T_1=\{
v\in T\mid r(s^{-1}(v))\cap C\ne \emptyset \}$, and for each $k> 1$, $T_k=\{ v\in T\setminus
T_{k-1}\mid r(s^{-1}(v))\cap T_{k-1}\ne \emptyset \}$.

Now, we will prove the result by induction on $k$ (the number of layers of $T$). For $k=1$, fix
$v\in T_1$, and notice that there exists $e\in s^{-1}(v)$ such that $r(e)\in C$. Let $\mu$ be a
cycle such that $r(e)\in \mu^0$, and let $f\in \mu^1$ such that $r(f)=r(e)$. By Condition (Sing)
there are no other edges satisfying the same hypotheses. Also, as $v\not\in C$, there are no edges
starting in $s(f)$ and ending in $v$. Thus, let $F^n_1=(E_n)[v \hookrightarrow s(f)]$. Since
$L(E_n)$ is purely infinite simple then $E_n$ necessarily satisfies Condition (L).  Thus Theorem
\ref{meu1_gen} applies to give $L(E_n)\cong L(F_1^n)$, and by the remark in the previous sentence,
$F_1^n$ satisfies Condition (Sing). By construction, $\vert E_n\vert = \vert F_1^n\vert$. Finally,
in $F_1^n$ we have $v\in C$, so that the set $T_1$ has one less element in $F_1^n$ than in $E_n$.
Applying this argument recurrently on the elements of the (finite) set $T_1$, we construct
$E_n^{(1)}$ satisfying Condition (Sing), with $L(E_n)\cong L(E_n^{(1)})$, and with the property
that the set of vertices in both graphs are identical, but the set $T$ in $E_n^{(1)}$ has one less
layer than the corresponding in $E_n$. Hence, the result holds by induction.

The parenthetical remark follows from the fact that Theorem \ref{meu1_gen} has been used in the
proof, and the isomorphism between
 Leavitt path algebras ensured by that result is not in general $\Z$-graded.

\end{proof}

Recall that for a ring $R$,  we denote by $K_0(R)$ the Grothendieck group of $R$.  This is the
group $F/S$, where $F$ is the free group generated by isomorphism classes of finitely generated
projective left $R$-modules, and $S$ is the subgroup of $F$ generated by symbols of the form
$[P\oplus Q]-[P]-[Q]$. As is standard, we denote the isomorphism class of $R$ in $K_0(R)$ by
$[1_R]$.  The group $K_0(R)$ is the universal group of the monoid $V(R)$ of isomorphism classes of
finitely generated projective left $R$-modules (with binary operation in $V(R)$ given by
$[A]+[B]=[A\oplus B]$). Because the rings we consider here are purely infinite simple Leavitt path
algebras, we have the following more explicit relationship between $K_0$ and $V$ in this setting:
$$V(L(E))\cong \{ 0\}\sqcup K_0(L(E))$$
(see for instance \cite[Corollary 2.2]{AGP}).

For a row-finite graph $E$, the {\it monoid of} $E$, denoted $M_E$, is the monoid generated by the
set $E^0$ of vertices of $E$ modulo appropriate relations, specifically,
$$M_E = \langle a_v, v\in E^0\mid a_v=\sum\limits
_{\{e\in s^{-1}(v)\}}a_{r(e)} \rangle.$$ It is shown in \cite[Theorem 2.5]{AMFP} that $V(L(E))\cong
M_E$ for any row-finite graph $E$. This yields $K_0(L(E))\cong \mbox{Grot}(M_E):=G$, where
$\mbox{Grot}(M_E)$ denotes the universal group of the monoid $M_E$. Since $M_E$ is finitely
generated, so is its universal group $G$. Thus $G$ admits a presentation $\pi:
\mathbb{Z}^n\rightarrow G$ (an epimorphism). Here $\mbox{ker}(\pi)$ is the subgroup of relations,
which in this setting corresponds to the image of the group homomorphism $A^t_E-I:
\mathbb{Z}^n\rightarrow\mathbb{Z}^n$, where $A^t_E$ is the transpose of the incidence matrix $A_E$
of $E$. Hence we get
$$K_0(L(E))\cong G\cong \mathbb{Z}^n
/\mbox{ker}(\pi)=\mathbb{Z}^n/\mbox{im}(A^t_E-I)=\mbox{coker}(A^t_E-I).$$ Moreover, under this
isomorphism the element $[1_{L(E)}]$ is represented by $(1,1,...,1)^t + \mbox{im}(A^t_E-I)$ in
$\mbox{coker}(A^t_E-I)$.

 Throughout the remainder of this
article we seek to describe properties of the Grothendieck groups $K_0(L(E))$ for various graphs
$E$.  To do so we will use the displayed isomorphism $K_0(L(E))\cong \mbox{coker}(A^t_E-I)$ often,
and without explicit mention.  (We present some examples below which indicate how one may directly
compute $\mbox{coker}(A^t_E-I)$.)

In the study of C$^*$-algebras, an important role is played by the Classification Theorem of purely
infinite simple unital nuclear C*-algebras (see e.g. \cite{Kirch, Phil}).  Specifically, Kirchberg
and Phillips (independently) showed that if $X$ and $Y$ are purely infinite simple unital
C$^*$-algebras (satisfying certain additional conditions), then $X\cong Y$ as C$^*$-algebras if and
only if (i) $K_0(X)\cong K_0(Y)$ via an isomorphism $\phi$ for which $\phi([1_X])=[1_Y]$, and (ii)
$K_1(X)\cong K_1(Y)$.

As it turns out, in the more specific case of purely infinite simple unital Cuntz-Krieger graph
C*-algebras, K-theoretic information is in fact encoded in the transpose $A_E^t$ of the incidence
matrix $A_E$ of the graph $E$. Specifically, when $E$ has no sinks, then by \cite[Theorem 3.9]{T1}
$$K_0(C^*(E))\cong \mbox{coker}(A^t_E-I) \  \hbox{  and } \ K_1(C^*(E))\cong \mbox{ker}(A^t_E-I),$$
where $I$ is the identity matrix of size $n=|E^0|$.

\medskip

We seek a similar result in the setting of purely infinite simple unital Leavitt path algebras. So
suppose $E$ and $F$ are finite graphs for which $L_K(E)$ and $L_K(F)$ are purely infinite simple
unital. By \cite[Theorem 11]{AA2} these graphs contain no sinks. By Proposition \ref{L:nosource} we
can assume without loss of generality that $E$ and $F$ have the same number $n$ of vertices and
that they have no sources. Thus if $K_0(L(E))\cong K_0(L(F))$, then using the previously
established isomorphism we get $\mbox{coker}(A^t_E-I) \cong \mbox{coker}(A^t_F-I)$. This in turn
implies (by the Fundamental Theorem of Finitely Generated Abelian Groups) the existence of
invertible matrices $P,Q \in \mbox{M}_n(\Z)$ such that $A^t_F-I=P(A^t_E-I)Q$. Thus
$\mbox{ker}(A^t_F-I)\cong \mbox{ker}(A^t_E-I)$ (as these are subgroups of $\mathbb{Z}^{n}$ having
equal rank); notice that in particular, since $K_1(C^*(E))\cong \mbox{ker}(A^t_E-I)$, we have
recovered the result of \cite[Theorem 3.9]{T1} for graph $C^*$-algebras. Moreover, by using the
unique unital ring map $\psi: \mathbb{Z}\rightarrow K$, we get that the $PAQ$-equivalence of
$A^t_E-I$ and $A^t_F-I$ also holds on $K$. If $K^{\times}$ denotes the multiplicative group on
nonzero elements in $K$, then the previous remark implies that $\mbox{coker}(A^t_E-I:
{(K^{\times})}^n \rightarrow {(K^{\times})}^n)$ and $\mbox{coker}(A^t_F-I: {(K^{\times})}^n
\rightarrow {(K^{\times})}^n)$ (where $A^t_E-I$ and $A^t_E-I$ are seen as multiplicative maps on
${(K^{\times})}^n$) are also isomorphic. Since by \cite[Theorem 3.19]{BrustengaPhD}, for any finite
graph $G$ with $n$ vertices with no sinks or sources we have
$$K_1(L(G))\cong \mbox{coker}(A^t_G-I: {(K^{\times})}^n
\rightarrow {(K^{\times})}^n)\oplus \mbox{ker}(A^t_G-I:\mathbb{Z}^n\rightarrow \mathbb{Z}^n),$$
 we conclude that the
hypothesis $K_0(L(E))\cong K_0(L(F))$ in fact yields $K_1(L(E))\cong K_1(L(F))$ as a consequence.
With this observation and the aforementioned Kirchberg - Phillips result in mind, it is then
natural to ask the following

\medskip
\bigskip

\noindent {\bf \textsc{The Classification Question for purely infinite simple unital Leavitt path
algebras.}}  Suppose $E$ and $F$ are graphs for which $L(E)$ and $L(F)$ are purely infinite simple
unital. If $K_0(L(E))\cong K_0(L(F))$ via an isomorphism $\phi$ having
$\phi([1_{L(E)}])=[1_{L(F)}]$, must $L(E)$ and $L(F)$ be isomorphic?

\medskip
\bigskip

Much of the remainder of this article is taken up in addressing The Classification Question.  We
notationally abbreviate the statement
 $$K_0(L(E))\cong K_0(L(F))\mbox{ via an isomorphism }\phi \mbox{ having }\phi([1_{L(E)}])=[1_{L(F)}]$$
 by writing
 $$(K_0(L(E)),[1_{L(E)}])\cong (K_0(L(F)),[1_{L(F)}]).$$

\smallskip

\begin{defi}\label{Leavittalgebra}
{\rm We recall that for each integer $n\geq 2$, the {\it Leavitt algebra} $L_n$ is the free
associative $K$-algebra with generators $\{x_i,y_i:1\leq i\leq n\}$ and relations
$$ (1)\  x_iy_j=\delta_{ij} \hbox{ for all } 1\leq i,j\leq n, \ \ \hbox{ and } \ \ (2)\ \sum_{i=1}^n y_ix_i=1 .$$}
\end{defi}
See \cite{AA1} or \cite{L} for additional information about $L_n$.   In particular, the isomorphism
$$L_n \cong L(R_n)$$
follows immediately, where $R_n$ is the ``rose with $n$ petals" graph

$$
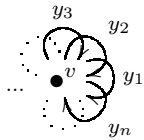
\xymatrix{ & {\bullet^v} \ar@(ur,dr) ^{y_1} \ar@(u,r) ^{y_2} \ar@(ul,ur) ^{y_3} \ar@{.} @(l,u) \ar@{.} @(dr,dl)
\ar@(r,d) ^{y_n} \ar@{}[l] ^{\ldots} }$$

\medskip

For $n\geq 2$ and $k\geq 1$ we define the graph $B_n^k$ to be
  $$
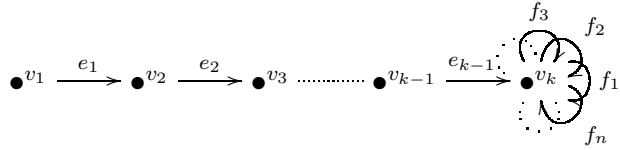
\xymatrix{{\bullet}^{v_1} \ar [r] ^{e_1} &
{\bullet}^{v_2} \ar [r] ^{e_2} & {\bullet}^{v_3} \ar@{.}[r] & {\bullet}^{v_{k-1}} \ar [r]
^{e_{k-1}} & {\bullet}^{v_k}
 \ar@(ur,dr) ^{f_1} \ar@(u,r) ^{f_2} \ar@(ul,ur) ^{f_3} \ar@{.} @(l,u) \ar@{.} @(dr,dl) \ar@(r,d) ^{f_n}& }$$
Then by \cite[Proposition 13]{AA2} we have
$$L(B_n^k)\cong {\rm M}_k(L_n).$$ We will use
this isomorphism throughout the sequel, often without explicit mention.

By \cite[Theorem 4.2]{AGP} we have that $K_0(L_n)\cong \Z/(n-1)\Z$. In fact, it is clear from this
isomorphism that $(K_0(L_n),[1_{L_n}])\cong (\Z/(n-1)\Z,\overline{1})$. Because $K_0$ is a Morita
invariant, we also necessarily have $K_0({\rm M}_k(L_n))\cong \Z/(n-1)\Z$ for any $k\in \N$. It is
straightforward to show that this isomorphism gives
$$(K_0({\rm M}_k(L_n)),[1_{{\rm M}_k(L_n)}])\cong (\Z/(n-1)\Z,\overline{k}).$$
   We will revisit this
isomorphism later, in two regards.  First, we will show in Example \ref{K0Lnfromgraphs} that it can
be re-established using tools from Leavitt path algebras.  Second, we will establish in Section
\ref{graphs} an affirmative answer to The Classification Question among a specific class of Leavitt
path algebras $L(E)$, to wit, if $(K_0(L(E)),[1_{L(E)}])\cong (\Z/(n-1)\Z,\overline{k})$, then
$L(E)\cong{\rm M}_k(L_n)$.

We now present some examples in which we explicitly compute $\mbox{coker}(A^t_E-I)$ for various
graphs $E$.  Additional examples and computations of this type can be found in \cite[page 32 and
Example 3.31]{T1}.

\begin{exem}\label{computeV}
{\rm Consider the graph
$${E_{1}^6}: \quad {
\def\labelstyle{\displaystyle}
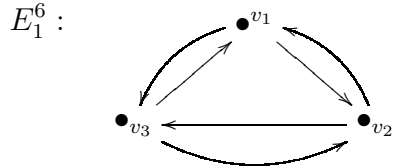
\xymatrix{ {} & \bullet^{v_1}  \ar[rd] \ar@/^{-10pt}/ [ld] &  {} \\
\bullet_{v_3}  \ar[ru] \ar@/^{-15pt}/ [rr]&  & \bullet_{v_2}
 \ar[ll]
\ar@/^{-10pt}/ [lu] \\
}}
$$
\vspace{.2truecm}

We compute $\mbox{coker}(A^t_{E_{1}^6}-I)$. First,
$$A^t_{E_{1}^6}-I=
\left(%
\begin{array}{ccc}
  -1 & 1 & 1 \\
  1 & -1 & 1 \\
  1 & 1 & -1 \\
\end{array}%
\right)
$$
Then, applying the classical $PAQ$-reduction, we get that $A^t_{E_{1}^6}-I$ is equivalent to the
diagonal matrix
$$D=
\left(%
\begin{array}{ccc}
  1 & 0 & 0 \\
  0 & 2 & 0 \\
  0 & 0 & 2 \\
\end{array}%
\right)
$$
while the invertible matrix $P$, which fixes the basis change in the arrival free group, is
$$P=
\left(%
\begin{array}{ccc}
  -1 & 0 & 0 \\
  1 & 0 & 1 \\
  1 & 1 & 0 \\
\end{array}%
\right)
$$
Then, $K_0(L({E_{1}^6}))\cong \mbox{coker}(D)=\mathbb{Z}/2\mathbb{Z}\oplus \mathbb{Z}/2\mathbb{Z}$.
On the other side, as in $\mbox{coker}(A^t_{E_{1}^6}-I)$ the element $[1_{L({E_{1}^6})}]$ is
represented by $(1,1,1)^t$, applying the change of basis we get that the image of the element
$[1_{L({E_{1}^6})}]$ in $\mbox{coker}(D)$ is $P\cdot (1,1,1)^t= (-1,2,2)^t$, modulo the relation
defined by $\mbox{im}(D)$, so that we conclude that $[1_{L({E_{1}^6})}]$ corresponds to
$(\overline{0}, \overline{0})$.}
\end{exem}

\begin{exem}\label{K0Lnfromgraphs}
{\rm Consider the graph $R_n$
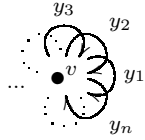
$$\xymatrix{ & {\bullet^v} \ar@(ur,dr) ^{y_1} \ar@(u,r) ^{y_2} \ar@(ul,ur) ^{y_3} \ar@{.} @(l,u) \ar@{.} @(dr,dl)
\ar@(r,d) ^{y_n} \ar@{}[l] ^{\ldots} }$$ and recall that $L(R_n)\cong L_n$. We will use the
$K_0$-picture described above to compute $K_0(L_n)$. We first compute $\mbox{coker}(A^t_{R_n}-I)$.
This is obvious, as
$$A^t_{R_n}-I=
\left(%
\begin{array}{c}
  n-1  \\
\end{array}%
\right)
$$
whence $K_0(L_n)\cong \mathbb{Z}/(n-1)\mathbb{Z}$. Since this matrix is in reduced form,
$[1_{L_n}]$ corresponds to $\overline{1}\in \mathbb{Z}/(n-1)\mathbb{Z}$.

Now consider  the graph $B_n^k$
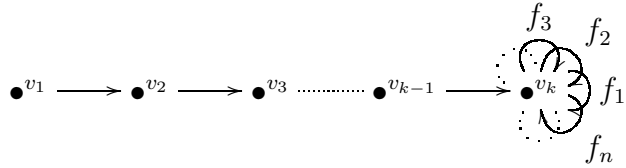
$$
{
\def\labelstyle{\displaystyle}
\xymatrix{{\bullet}^{v_1} \ar [r]  & {\bullet}^{v_2} \ar [r] & {\bullet}^{v_3} \ar@{.}[r] &
{\bullet}^{v_{k-1}} \ar [r] & {\bullet}^{v_k} \ar@(ur,dr) ^{f_1} \ar@(u,r) ^{f_2} \ar@(ul,ur)
^{f_3} \ar@{.} @(l,u) \ar@{.} @(dr,dl) \ar@(r,d) ^{f_n}& }}
$$
and recall that $L(B_n^k)\cong {\rm M}_k(L_n)$. First,
$$A^t_{B_n^k}-I=
\left(%
\begin{array}{cccccc}
  n-1 & 1 & 0 & 0 & \cdots & 0 \\
  0 & -1 & 1 & 0 & \cdots & 0 \\
  0 & 0 & -1 & 1 & \cdots & 0 \\
  \vdots & \vdots  & \vdots & \ddots & \cdots & \vdots \\
  0 & 0 & 0 & \cdots & -1& 1 \\
  0 & 0 & 0 & \cdots & 0 & -1 \\
\end{array}%
\right)
$$
Then, applying the classical $PAQ$-reduction, we get that $A^t_{B_n^k}-I$ is equivalent to the
diagonal matrix
$$D=
\left(%
\begin{array}{cccccc}
  n-1 & 0 & 0 & 0 & \cdots & 0 \\
  0 & 1 & 0 & 0 & \cdots & 0 \\
  0 & 0 & 1 & 0 & \cdots & 0 \\
  \vdots & \vdots  & \vdots & \ddots & \cdots & \vdots \\
  0 & 0 & 0 & \cdots & 1& 0 \\
  0 & 0 & 0 & \cdots & 0 & 1 \\
\end{array}%
\right)
$$
while the invertible matrix $P$, which fixes the basis change in the arrival free group, is
$$P=
\left(%
\begin{array}{cccccc}
  1 & 1 & 1 & 1 & \cdots & 1 \\
  0 & -1 & -1 & -1 & \cdots & -1 \\
  0 & 0 & -1 & -1 & \cdots & -1 \\
  \vdots & \vdots  & \vdots & \ddots & \cdots & \vdots \\
  0 & 0 & 0 & \cdots & -1& -1 \\
  0 & 0 & 0 & \cdots & 0 & -1 \\
\end{array}%
\right)
$$
Then $K_0(L(B_n^k))\cong \mbox{coker}(D)=\mathbb{Z}/(n-1)\mathbb{Z}$. (We have thereby
re-established a
 previously observed isomorphism between $K_0(L(B_n^k))$ and
$\mathbb{Z}/(n-1)\mathbb{Z}$.)  On the other side, as in $\mbox{coker}(A^t_{B_n^k}-I)$ the element
$[1_{B_n^k}]$ is represented by $(1,1,\dots ,1)^t$, applying the change of basis we get that the
image of the element $[1_{L({B_n^k})}]$ in $\mbox{coker}(D)$ is $P\cdot (1,1,\dots ,1)^t=
(k,-(k-1), -(k-2), \dots,-2, -1)^t$, modulo the relation defined by $\mbox{im}(D)$, so that we
conclude that $[1_{L({B_n^k})}]$ corresponds to $\overline{k}$. }
\end{exem}

\section{Graphs}\label{graphs}

In this section we will show how to use Theorems \ref{meu1_gen} and \ref{outsplitiso} in order to
classify purely infinite simple unital Leavitt path algebras according to their $K_0$-data.
Specifically, we give an affirmative answer to The Classification Question for Leavitt path
algebras coming from various collections of graphs.

  Among these collections will be graphs whose
$K_0$-data matches the $K_0$-data for Leavitt path algebras of the form $L(E)\cong {\rm M}_k(L_n)$.
(For instance, as noted previously, by \cite[Proposition 13]{AA2} the graph $B_n^k$ has
$L(B_n^k)\cong  {\rm M}_k(L_n)$.)  That is, we will have Leavitt path algebras $L(E)$ for which
$(K_0(L(E)), [1_{L(E)}])\cong (\mathbb{Z}/(n-1)\mathbb{Z},\overline{k})$. For such collections we
provide additional evidence that The Classification Question has an affirmative answer, by showing
that the relevant Leavitt path algebras are indeed isomorphic to ${\rm M}_k(L_n)$ for appropriate
$n,k$.  This is shown at the end of each of the germane subsections into which the section is
divided.

We note that all of the graphs we consider throughout this section satisfy Condition (L) (since the
graphs arise in the context of purely infinite simple Leavitt path algebras).  Thus Theorem
\ref{meu1_gen} applies to all of the shift graph constructions produced here, regardless of the
size of the field $K$.

\subsection{Graphs with two vertices}\label{ssect_2vert}

We start by analyzing graphs having two vertices, which satisfy Condition  (Sing), and for which
the associated Leavitt path algebra is purely infinite simple.  Concretely, they are the following:

$$
\widehat{R_2}: {
\def\labelstyle{\displaystyle}
\xymatrix{\bullet^{v_1}\dloopd{} \ar@/^1pc/ [r] & {\bullet}^{v_2} \ar@/^1pc/ [l] \uloopd{} }};
\qquad {
\def\labelstyle{\displaystyle}S_2:\quad \xymatrix{
{\bullet}^{v_1} \ar@(ul,dl)  \ar@/^1pc/ [r] & {\bullet}^{v_2} \ar@/^1pc/ [l]  }}; \qquad
{B_2^2}:\quad {
\def\labelstyle{\displaystyle}
\xymatrix{ \bullet^{v_1}\uloopr{}\dloopr{} & \bullet^{v_2} \ar[l]}}.
$$
Using the description of $(K_0(L(E)), [1_{L(E)}])$ given at the end of Section
\ref{purelyinfinitesimple}, it is straightforward to show that each of these three graphs has
$(K_0(L(E)), [1_{L(E)}])\cong (\{ \overline{0}\}, \overline{0})$.   But then the isomorphisms
between the respective algebras can be found in Examples \ref{isoofshiftexample} and
\ref{outsplitexample}.

Thus we have answered in the affirmative a specific case of The Classification Question for purely
infinite simple unital Leavitt path algebras.

\begin{prop}\label{twovertices}
Suppose $E$ and $F$ are graphs having Condition (Sing), for which $L(E)$ and $L(F)$ are purely
infinite simple unital, and $|E^0| = |F^0| = 2$. If $K_0(L(E))\cong K_0(L(F))$ via an isomorphism
$\phi$ for which $\phi([1_{L(E)}])=[1_{L(F)}]$, then $L(E)\cong L(F)$.
\end{prop}

The three graphs of Proposition \ref{twovertices} each have $K_0$-data $(\{ \overline{0}\},
\overline{0})$, which matches the $K_0$-data of ${\rm M}_1(L_2)\cong L_2$. We show that in fact
$L(E)\cong L_2$ for each of these three graphs. This will follow directly from the isomorphism
$L_2\cong L(B_2^2)$ ensured by Propositions \ref{L:congruent mod n-1} and \ref{P:dividingpower}
below.

\subsection{Graphs with three vertices}

We continue by analyzing graphs having three vertices, which satisfy Condition  (Sing), and for
which the associated Leavitt path algebra is purely infinite simple.   It turns out there exist
$34$ such graphs.  Unlike the previously analyzed situation for graphs with two vertices, there
will be more than one pair of the form $(K_0(L(E)),[1_{L(E)}])$ arising from this collection.
(There are seven such pairs, to be exact.)  We partition all 34 of these graphs along the seven
$K_0$-data pairs, and then use the tools of Section \ref{XXX} to show that the Leavitt path
algebras within each equivalence class are indeed pairwise isomorphic.  Throughout we use without
mention the description of $(K_0(L(E)),[1_{L(E)}])$ presented at the end of Section
\ref{purelyinfinitesimple}.

\bigskip

\begin{point} {\rm \textbf{ $(K_0(L(E)), [1_{L(E)}])\cong
(\{\overline{0}\},\overline{0})$ }: In this situation we have $18$ graphs, listed as follows:
$$ {E_1^1}: \quad {
\def\labelstyle{\displaystyle}
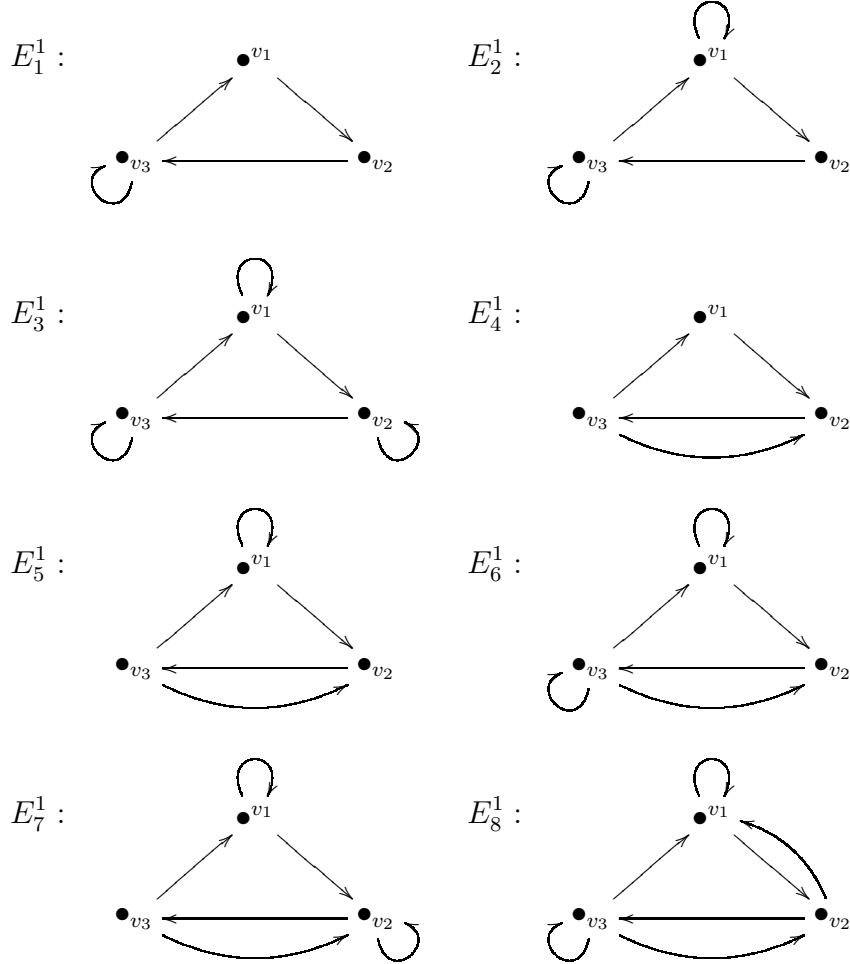
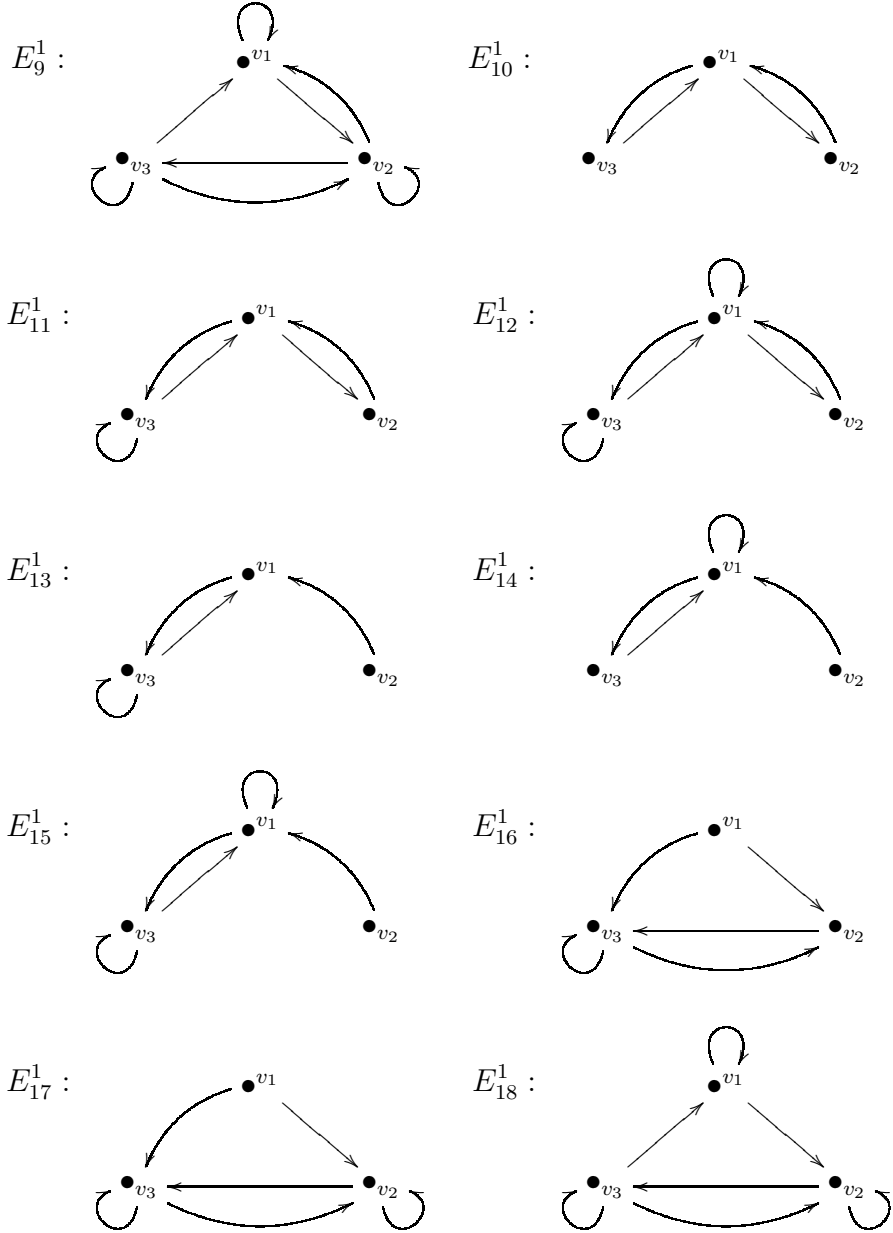
\xymatrix{ {} & \bullet^{v_1}  \ar[rd]  &  {} \\
\bullet_{v_3} \ar@(d,l) \ar[ru] &  & \bullet_{v_2}  \ar[ll]
 \\
}} \qquad {E_2^1}: \quad {
\def\labelstyle{\displaystyle}
\xymatrix{ {} & \bullet^{v_1} \uloopr{} \ar[rd]  &  {} \\
\bullet_{v_3} \ar@(d,l) \ar[ru] &  & \bullet_{v_2} \ar[ll]
\\
}}
$$
\vspace{.1truecm}

$$ {E_3^1}: \quad {
\def\labelstyle{\displaystyle}
\xymatrix{ {} & \bullet^{v_1} \uloopr{} \ar[rd]  &  {} \\
\bullet_{v_3} \ar@(d,l) \ar[ru] &  & \bullet_{v_2} \ar@(d,r) \ar[ll]
\\
}} \qquad {E_4^1}: \quad {
\def\labelstyle{\displaystyle}
\xymatrix{ {} & \bullet^{v_1}  \ar[rd]  &  {} \\
\bullet_{v_3}  \ar[ru] \ar@/^{-15pt}/ [rr]&  & \bullet_{v_2}
 \ar[ll]
 \\
}}
$$
\vspace{.1truecm}

$$ {E_5^1}: \quad {
\def\labelstyle{\displaystyle}
\xymatrix{ {} & \bullet^{v_1} \uloopr{} \ar[rd] &  {} \\
\bullet_{v_3} \ar[ru] \ar@/^{-15pt}/ [rr]&  & \bullet_{v_2}  \ar[ll]
\\
}} \qquad {E_6^1}: \quad {
\def\labelstyle{\displaystyle}
\xymatrix{ {} & \bullet^{v_1} \uloopr{} \ar[rd] &  {} \\
\bullet_{v_3} \ar@(d,l) \ar[ru] \ar@/^{-15pt}/ [rr]&  & \bullet_{v_2} \ar[ll]
\\
}}
$$
\vspace{.1truecm}

$$ {E_7^1}: \quad {
\def\labelstyle{\displaystyle}
\xymatrix{ {} & \bullet^{v_1} \uloopr{} \ar[rd] &  {} \\
\bullet_{v_3} \ar[ru] \ar@/^{-15pt}/ [rr]&  & \bullet_{v_2} \ar@(d,r) \ar[ll]
\\
}} \qquad {E_8^1}: \quad {
\def\labelstyle{\displaystyle}
\xymatrix{ {} & \bullet^{v_1} \uloopr{} \ar[rd] &  {} \\
\bullet_{v_3} \ar@(d,l) \ar[ru] \ar@/^{-15pt}/ [rr]&  & \bullet_{v_2}  \ar[ll]
\ar@/^{-10pt}/ [lu] \\
}}
$$
\vspace{.1truecm}

$$ {E_9^1}: \quad {
\def\labelstyle{\displaystyle}
\xymatrix{ {} & \bullet^{v_1} \uloopr{} \ar[rd]  &  {} \\
\bullet_{v_3} \ar@(d,l) \ar[ru] \ar@/^{-15pt}/ [rr]&  & \bullet_{v_2} \ar@(d,r) \ar[ll]
\ar@/^{-10pt}/ [lu] \\
}} \qquad {E_{10}^1}: \quad {
\def\labelstyle{\displaystyle}
\xymatrix{ {} & \bullet^{v_1} \ar[rd] \ar@/^{-10pt}/ [ld] &  {} \\
\bullet_{v_3} \ar[ru] &  & \bullet_{v_2}
\ar@/^{-10pt}/ [lu] \\
}}
$$
\vspace{.1truecm}

$$ {E_{11}^1}: \quad {
\def\labelstyle{\displaystyle}
\xymatrix{ {} & \bullet^{v_1}  \ar[rd] \ar@/^{-10pt}/ [ld] &  {} \\
\bullet_{v_3} \ar@(d,l) \ar[ru] &  & \bullet_{v_2}
\ar@/^{-10pt}/ [lu] \\
}} \qquad {E_{12}^1}: \quad {
\def\labelstyle{\displaystyle}
\xymatrix{ {} & \bullet^{v_1} \uloopr{} \ar[rd] \ar@/^{-10pt}/ [ld] &  {} \\
\bullet_{v_3} \ar@(d,l) \ar[ru] &  & \bullet_{v_2}
\ar@/^{-10pt}/ [lu] \\
}}
$$
\vspace{.1truecm}

$$ {E_{13}^1}: \quad {
\def\labelstyle{\displaystyle}
\xymatrix{ {} & \bullet^{v_1}  \ar@/^{-10pt}/ [ld] &  {} \\
\bullet_{v_3} \ar@(d,l) \ar[ru] &  & \bullet_{v_2}
\ar@/^{-10pt}/ [lu] \\
}} \qquad {E_{14}^1}: \quad {
\def\labelstyle{\displaystyle}
\xymatrix{ {} & \bullet^{v_1} \uloopr{}  \ar@/^{-10pt}/ [ld] &  {} \\
\bullet_{v_3}  \ar[ru] &  & \bullet_{v_2}
\ar@/^{-10pt}/ [lu] \\
}}
$$
\vspace{.1truecm}

$$ {E_{15}^1}: \quad {
\def\labelstyle{\displaystyle}
\xymatrix{ {} & \bullet^{v_1} \uloopr{} \ar@/^{-10pt}/ [ld] &  {} \\
\bullet_{v_3} \ar@(d,l) \ar[ru] &  & \bullet_{v_2}
\ar@/^{-10pt}/ [lu] \\
}} \qquad {E_{16}^1}: \quad {
\def\labelstyle{\displaystyle}
\xymatrix{ {} & \bullet^{v_1}  \ar[rd] \ar@/^{-10pt}/ [ld] &  {} \\
\bullet_{v_3} \ar@(d,l) \ar@/^{-15pt}/ [rr]&  & \bullet_{v_2} \ar[ll]
 \\
}}
$$
\vspace{.1truecm}

$$ {E_{17}^1}: \quad {
\def\labelstyle{\displaystyle}
\xymatrix{ {} & \bullet^{v_1} \ar[rd] \ar@/^{-10pt}/ [ld] &  {} \\
\bullet_{v_3} \ar@(d,l)  \ar@/^{-15pt}/ [rr]&  & \bullet_{v_2} \ar@(d,r) \ar[ll]
 \\
}} \qquad {E_{18}^1}: \quad {
\def\labelstyle{\displaystyle}
\xymatrix{ {} & \bullet^{v_1} \uloopr{} \ar[rd]  &  {} \\
\bullet_{v_3} \ar@(d,l) \ar[ru] \ar@/^{-15pt}/ [rr]&  & \bullet_{v_2} \ar@(d,r) \ar[ll]
\\
}}
$$
\vspace{.2truecm}

Now, we prove the isomorphisms as follows: First consider the out-splitting of
$$
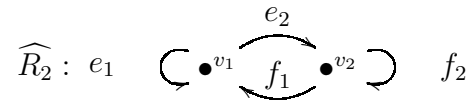
\widehat{R_2}: {
\def\labelstyle{\displaystyle}
\xymatrix{{e_1} &\bullet^{v_1}\dloopd{} \ar@/^1pc/ [r]^{e_2} & {\bullet}^{v_2} \ar@/^1pc/ [l]_{f_1}
\uloopd{} &{f_2} }}$$ partitioning the edges in $\mathcal{P}=\{ e_1\}\cup\{ e_2\}\cup\{ f_1,
f_2\}$, and notice that $(\widehat{R_2})_s(\mathcal{P})=E_6^1$. Thus, $L(\widehat{R_2})\cong
L(E_6^1)$ by Theorem \ref{outsplitiso}. Now, consider the out-splitting of
$${
\def\labelstyle{\displaystyle}S_2:\quad 
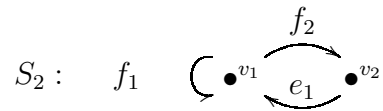
\xymatrix{{f_1}&
{\bullet}^{v_1} \ar@(ul,dl)  \ar@/^1pc/ [r]^{f_2} & {\bullet}^{v_2} \ar@/^1pc/ [l]_{e_1}  }}$$
partitioning the edges in $\mathcal{P}=\{ e_1\}\cup\{ f_1\}\cup\{ f_2\}$, and notice that
$(S_2)_s(\mathcal{P})=E_5^1$. Thus, $L({S_2})\cong L(E_5^1)$ by Theorem \ref{outsplitiso}. Finally,
the maximal out-splitting of $B_2^2$ equals $E_{17}^1$, whence $L(B_2^2)\cong L(\widehat{B_2^2}) =
L(E_{17}^1)$  by Corollary  \ref{isomaxoutsplit}. Since $L(S_2)\cong L(\widehat{R_2})\cong
L(B_2^2)$ by Proposition \ref{twovertices}, we have shown that $L(E_5^1)\cong L(E_6^1) \cong
L(E_{17}^1)$, which in turn can be used to verify the isomorphisms with all the remaining indicated
Leavitt path algebras by noticing that

\begin{enumerate}
\item $E_{17}^1(v_3\hookrightarrow v_2)=E_{16}^1$
\item $E_{16}^1(v_3\hookrightarrow v_1)=E_4^1$ \item
$E_{16}^1(v_1\hookrightarrow v_3)=E_{14}^1$ \item $E_{2}^1(v_2\hookrightarrow v_3)=E_{5}^1$ \item
$E_{1}^1(v_2\hookrightarrow v_3)=E_{4}^1$ \item $E_{5}^1(v_3\hookrightarrow v_1)=E_{10}^1$
\item $E_{7}^1(v_1\hookrightarrow v_3)=E_{2}^1$ \item
$E_{7}^1(v_3\hookrightarrow v_1)=E_{11}^1$ \item $E_{9}^1(v_3\hookrightarrow v_2)=E_{6}^1$ \item
$E_{9}^1(v_2\hookrightarrow v_3)=E_{12}^1$ \item $E_{8}^1(v_1\hookrightarrow v_3)=E_{7}^1$ \item
$E_{14}^1(v_3\hookrightarrow v_2)=E_{13}^1$ \item $E_{15}^1(v_1\hookrightarrow v_3)=E_{14}^1$ \item
$E_{18}^1(v_2\hookrightarrow v_3)=E_{7}^1$ \item $E_{18}^1(v_1\hookrightarrow v_3)=E_{3}^1$
\end{enumerate}
 Then, all those Leavitt path algebras are pairwise isomorphic
by Theorem \ref{meu1_gen}, so we are done.

\medskip

The eighteen graphs of this subsection have $K_0$-data $(\{\overline{0}\},\overline{0})$. But the
purely infinite simple Leavitt path algebra $L(B_2^2)\cong L_2$ has this same $K_0$-data as well.
As further evidence of an affirmative answer to The Classification Question, we note that indeed we
have shown, for all eighteen graphs $E$ in this subsection, that $L(E)\cong L_2$.  (We established
the isomorphism $L(B_2^2)\cong L(E_{17}^1)$ in the course of the proof.)

}
\end{point}

\bigskip

\begin{point} {\rm \textbf{$(K_0(L(E)), [1_{L(E)}])\cong (\mathbb{Z}/2\mathbb{Z},
\overline{0})$}: In this situation we have $6$ graphs, listed as follows:
$${E_{1}^2}: \quad {
\def\labelstyle{\displaystyle}
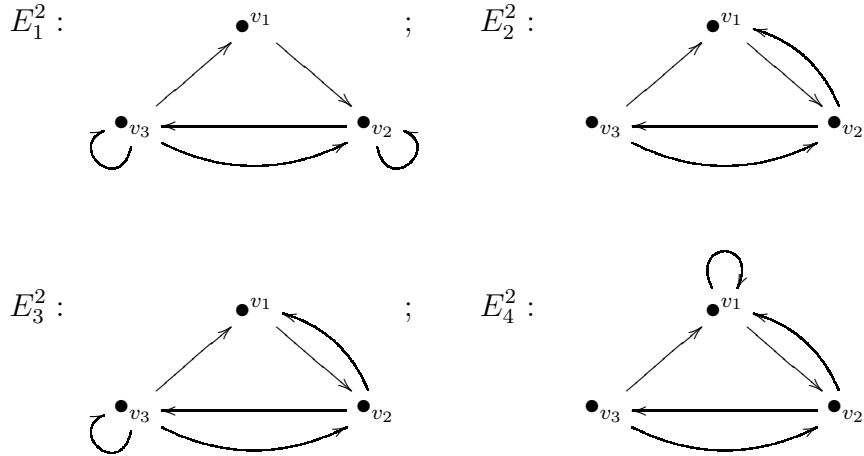
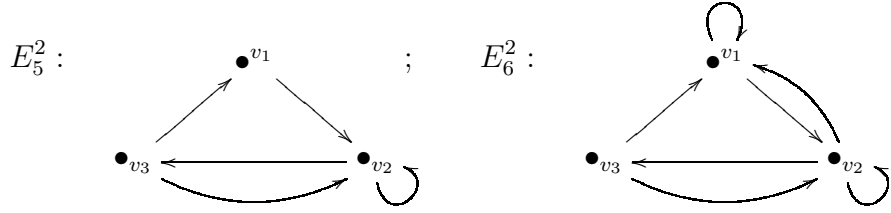
\xymatrix{ {} & \bullet^{v_1} \ar[rd]  &  {} \\
\bullet_{v_3} \ar@(d,l) \ar[ru] \ar@/^{-15pt}/ [rr]&  & \bullet_{v_2} \ar@(d,r) \ar[ll]
 \\
}}; \qquad {E_{2}^2}: \quad {
\def\labelstyle{\displaystyle}
\xymatrix{ {} & \bullet^{v_1} \ar[rd]  &  {} \\
\bullet_{v_3}  \ar[ru] \ar@/^{-15pt}/ [rr]&  & \bullet_{v_2}
 \ar[ll]
\ar@/^{-10pt}/ [lu] \\
}}
$$
\vspace{.1truecm}

$${E_{3}^2}: \quad {
\def\labelstyle{\displaystyle}
\xymatrix{ {} & \bullet^{v_1}  \ar[rd]  &  {} \\
\bullet_{v_3} \ar@(d,l) \ar[ru] \ar@/^{-15pt}/ [rr]&  & \bullet_{v_2}  \ar[ll]
\ar@/^{-10pt}/ [lu] \\
}}; \qquad {E_{4}^2}: \quad {
\def\labelstyle{\displaystyle}
\xymatrix{ {} & \bullet^{v_1} \uloopr{} \ar[rd] &  {} \\
\bullet_{v_3}  \ar[ru] \ar@/^{-15pt}/ [rr]&  & \bullet_{v_2}
 \ar[ll]
\ar@/^{-10pt}/ [lu] \\
}}
$$
\vspace{.1truecm}

$${E_{5}^2}: \quad {
\def\labelstyle{\displaystyle}
\xymatrix{ {} & \bullet^{v_1}  \ar[rd]  &  {} \\
\bullet_{v_3}  \ar[ru] \ar@/^{-15pt}/ [rr]&  & \bullet_{v_2} \ar@(d,r) \ar[ll]
 \\
}}; \qquad {E_{6}^2}: \quad {
\def\labelstyle{\displaystyle}
\xymatrix{ {} & \bullet^{v_1} \uloopr{} \ar[rd]  &  {} \\
\bullet_{v_3}  \ar[ru] \ar@/^{-15pt}/ [rr]&  & \bullet_{v_2} \ar@(d,r) \ar[ll]
\ar@/^{-10pt}/ [lu] \\
}}
$$
\vspace{.2truecm}

Now, we prove the isomorphisms as follows: First notice that
\begin{enumerate}
\item $E_{6}^2(v_3\hookrightarrow v_1)=E_{3}^2$ \item
$E_{6}^2(v_1\hookrightarrow v_3)=E_{1}^2$ \item $E_{1}^2(v_2\hookrightarrow v_3)=E_{5}^2$ \item
$E_{4}^2(v_1\hookrightarrow v_3)=E_{5}^2$ \item $E_{5}^2(v_1\hookrightarrow v_2)=E_{2}^2$
\end{enumerate}
Then, all those Leavitt path algebras are pairwise isomorphic by Theorem \ref{meu1_gen}, so we are
done.

\medskip

The six graphs of this subsection have $K_0$-data $(\Z/2\Z,\overline{0}) = (\Z/2\Z,\overline{2})$.
But the purely infinite simple Leavitt path algebra ${\rm M}_2(L_3)$ has this same $K_0$-data as
well.  As further evidence of an affirmative answer to The Classification Question, we now show,
for all six graphs $E$ in this subsection, that $L(E)\cong {\rm M}_2(L_3)$.   To see this, by
\cite[Proposition 13]{AA2} we have ${\rm M}_2(L_3)\cong L(B_3^2)$, where
$${B_3^2}:\quad {
\def\labelstyle{\displaystyle}
\xymatrix{ \bullet^{v_1}\ar@(ul,dl)_{(3)} & \bullet^{v_2} \ar[l]}}$$ (here the notation $(n)$
indicates that there are $n$ parallel edges), and a single application of Theorem \ref{meu1_gen}
gives us $E_1={B_3^2}(v_1\hookrightarrow v_2)$, where
$${
\def\labelstyle{\displaystyle}E_1:\quad \xymatrix{
(2) & &{\bullet}^{v_1} \ar@(ul,dl)_{e_2,e_3} \ar@/^1pc/ [r]^{e_1} & {\bullet}^{v_2} \ar@/^1pc/
[l]_{f}  }}.$$ Partitioning the edges in $\mathcal{P}=\{ f\}\cup\{ e_1, e_2\}\cup\{ e_3\}$, we get
$(E_1)_s(\mathcal{P})=E_6^2$, so that the result holds by Theorem \ref{outsplitiso}, as desired.}
\end{point}

\bigskip

\begin{point}{\rm \textbf{$(K_0(L(E)), [1_{L(E)}])\cong (\mathbb{Z}/2\mathbb{Z},
\overline{1})$}: In this situation we have $4$ graphs, listed as follows:
$${E_{1}^3}: \quad {
\def\labelstyle{\displaystyle}
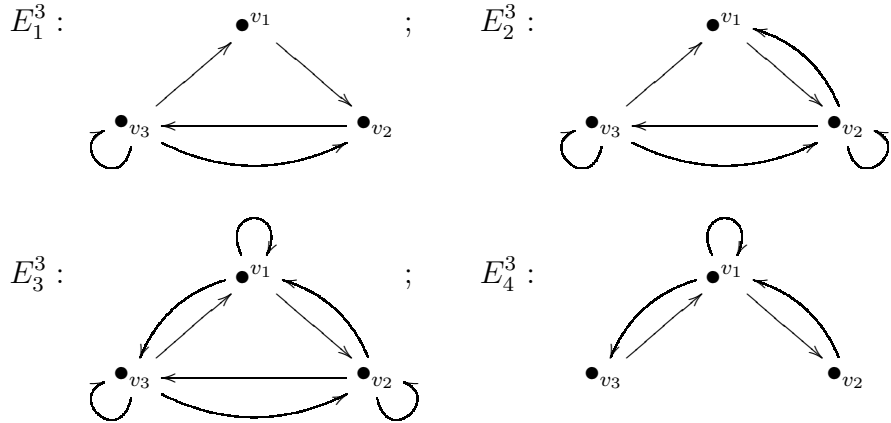
\xymatrix{ {} & \bullet^{v_1}  \ar[rd]  &  {} \\
\bullet_{v_3} \ar@(d,l) \ar[ru] \ar@/^{-15pt}/ [rr]&  & \bullet_{v_2} \ar[ll]
\\
}}; \qquad {E_{2}^3}: \quad {
\def\labelstyle{\displaystyle}
\xymatrix{ {} & \bullet^{v_1}  \ar[rd]  &  {} \\
\bullet_{v_3} \ar@(d,l) \ar[ru] \ar@/^{-15pt}/ [rr]&  & \bullet_{v_2} \ar@(d,r) \ar[ll]
\ar@/^{-10pt}/ [lu] \\
}}
$$
\vspace{.1truecm}

$${E_{3}^3}: \quad {
\def\labelstyle{\displaystyle}
\xymatrix{ {} & \bullet^{v_1} \uloopr{} \ar[rd] \ar@/^{-10pt}/ [ld] &  {} \\
\bullet_{v_3} \ar@(d,l) \ar[ru] \ar@/^{-15pt}/ [rr]&  & \bullet_{v_2} \ar@(d,r) \ar[ll]
\ar@/^{-10pt}/ [lu] \\
}}; \qquad {E_{4}^3}: \quad {
\def\labelstyle{\displaystyle}
\xymatrix{ {} & \bullet^{v_1} \uloopr{} \ar[rd] \ar@/^{-10pt}/ [ld] &  {} \\
\bullet_{v_3}  \ar[ru] &  & \bullet_{v_2}
\ar@/^{-10pt}/ [lu] \\
}}
$$
\vspace{.2truecm}

Now, we prove the isomorphisms as follows: First notice that
\begin{enumerate}
\item $E_{3}^3(v_2\hookrightarrow v_1)=E_{2}^3$ \item
$E_{2}^3(v_2\hookrightarrow v_3)=E_{4}^3$ \item $E_{4}^3(v_3\hookrightarrow v_2)=E_{1}^3$
\end{enumerate}
Then, all those Leavitt path algebras are pairwise isomorphic by Theorem \ref{meu1_gen}, so we are
done.

\medskip

The four graphs of this subsection have $K_0$-data $(\Z/2\Z,\overline{1})$.  But the purely
infinite simple Leavitt path algebra ${\rm M}_1(L_3)\cong L_3$
 has this same $K_0$-data as well.  As further evidence of an affirmative answer to The Classification Question, we now show,
  for all four graphs $E$ in this subsection,
that $L(E)\cong L_3$.   To see this, recall that $L_3\cong L(R_3)$, where
$${R_3}:\quad {
\def\labelstyle{\displaystyle}
\xymatrix{\bullet^{v_1}\ar@(ul,dl)_{(3)}}}.$$ Notice that the maximal out-splitting of $R_3$ equals
$E_3^3$, whence by Theorem \ref{outsplitiso} we get the desired result.}
\end{point}

\bigskip

\begin{point}{\rm \textbf{$(K_0(L(E)), [1_{L(E)}])\cong (\mathbb{Z}/3\mathbb{Z},
\overline{1})$}: In this situation we have $2$ graphs, listed as follows:
$${E_{1}^4}: \quad {
\def\labelstyle{\displaystyle}
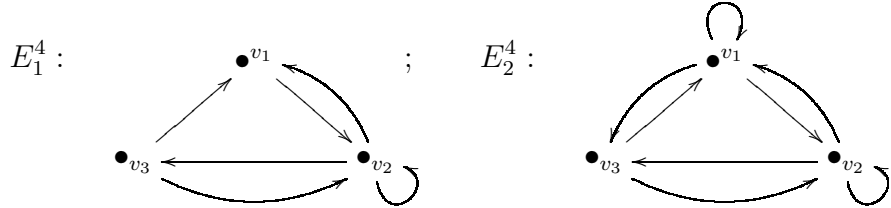
\xymatrix{ {} & \bullet^{v_1}  \ar[rd] &  {} \\
\bullet_{v_3}  \ar[ru] \ar@/^{-15pt}/ [rr]&  & \bullet_{v_2} \ar@(d,r) \ar[ll]
\ar@/^{-10pt}/ [lu] \\
}}; \qquad {E_{2}^4}: \quad {
\def\labelstyle{\displaystyle}
\xymatrix{ {} & \bullet^{v_1} \uloopr{} \ar[rd] \ar@/^{-10pt}/ [ld] &  {} \\
\bullet_{v_3}  \ar[ru] \ar@/^{-15pt}/ [rr]&  & \bullet_{v_2} \ar@(d,r) \ar[ll]
\ar@/^{-10pt}/ [lu] \\
}}
$$
\vspace{.2truecm}

By noticing that $E_{2}^4(v_2\hookrightarrow v_1)=E_{1}^4$, these Leavitt path algebras are
isomorphic by Theorem \ref{meu1_gen}, so we are done.

The two graphs of this subsection have $K_0$-data $(\Z/3\Z,\overline{1})$.  But the purely infinite
simple Leavitt path algebra ${\rm M}_1(L_4)\cong L_4$
 has this same $K_0$-data as well.  As further evidence of an affirmative answer to The Classification Question, we now show,
  for both graphs $E$ in this subsection,
that $L(E)\cong L_4$.    For, recall that $L_4\cong L(R_4)$, where

$${R_4}:\quad {
\def\labelstyle{\displaystyle}
\xymatrix{\bullet^{v_1}\ar@(ul,dl)_{(4)}}}.$$ Consider the maximal out-splitting $\widehat{R_4}$ of
$R_4$\vspace{.2truecm}

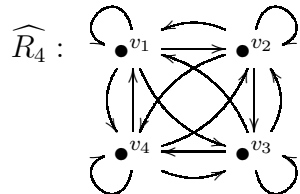
$${
\def\labelstyle{\displaystyle} \widehat{R_4}:\quad
\xymatrix{\bullet^{v_1} \ar@(u,l) \ar[r] \ar@/^{-10pt}/[rd] \ar@/^{-10pt}/ [d] & \bullet^{v_2}
\ar@(u,r) \ar@/^{-10pt}/ [l]
\ar@/^{-10pt}/[ld] \ar[d]\\
\bullet^{v_4} \ar@(d,l) \ar@/^{-10pt}/[ru] \ar[u] \ar@/^{-10pt}/ [r] & \bullet^{v_3} \ar@(d,r)
\ar[l] \ar@/^{-10pt}/ [u]
\ar@/^{-10pt}/[lu] \\
}}$$ \vspace{.1truecm}

Then, $L(R_4)\cong L(\widehat{R_4})$ by Corollary \ref{isomaxoutsplit}. Consider the graph
$\widehat{E}=[[\widehat{R_4}(v_1\hookrightarrow v_2)](v_1\hookrightarrow v_3)](v_1\hookrightarrow
v_4)$, that is,
$${
\def\labelstyle{\displaystyle} \widehat{E}:\quad
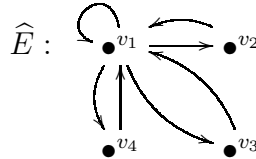
\xymatrix{\bullet^{v_1} \ar@(u,l) \ar[r] \ar@/^{-10pt}/[rd] \ar@/^{-10pt}/ [d] & \bullet^{v_2}
\ar@/^{-10pt}/ [l]
\\
\bullet^{v_4} \ar[u] & \bullet^{v_3} \ar@/^{-10pt}/[lu] \\
}}$$ \vspace{.1truecm}

Hence, $L(\widehat{R_4})\cong L(\widehat{E})$ by several applications of Theorem \ref{meu1_gen}. If
we separate the edges in $E_1^4$ in such way that the edges emitted by $v_3$ are divided in two
singletons, then by defining $F=(E_1^4)_s(\mathcal{P})$, we get \vspace{.2truecm}

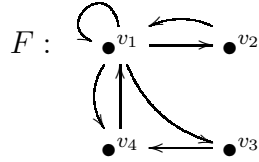
$${
\def\labelstyle{\displaystyle} {F}:\quad
\xymatrix{\bullet^{v_1} \ar@(u,l) \ar[r] \ar@/^{-10pt}/[rd] \ar@/^{-10pt}/ [d] & \bullet^{v_2}
\ar@/^{-10pt}/ [l]
\\
\bullet^{v_4} \ar[u] & \bullet^{v_3} \ar[l] \\
}}$$ \vspace{.1truecm}

Clearly, $F=\widehat{E}(v_3\hookrightarrow v_4)$. Thus, the result holds by Theorem
\ref{outsplitiso}.}
\end{point}

\bigskip

\begin{point}{\rm \textbf{$(K_0(L(E)), [1_{L(E)}])\cong (\mathbb{Z}/4\mathbb{Z},
\overline{2})$}: In this situation we have one graph, listed as follows:
$${E_{1}^5}: \quad {
\def\labelstyle{\displaystyle}
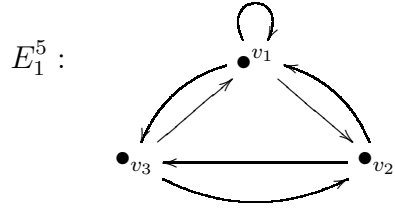
\xymatrix{ {} & \bullet^{v_1} \uloopr{} \ar[rd] \ar@/^{-10pt}/ [ld] &  {} \\
\bullet_{v_3}  \ar[ru] \ar@/^{-15pt}/ [rr]&  & \bullet_{v_2} \ar[ll]
\ar@/^{-10pt}/ [lu] \\
}}
$$
\vspace{.2truecm}

The one graph of this subsection has $K_0$-data $(\Z/4\Z,\overline{2})$.  But the purely infinite
simple Leavitt path algebra ${\rm M}_2(L_5)$
 has this same $K_0$-data as well.  As further evidence of an affirmative answer to The Classification Question, we now show,
  for the graph $E$ in this subsection,
that $L(E)\cong {\rm M}_2(L_5)$.
 For, let $E_1=E_1^5(v_2\hookrightarrow v_1)$,
$${
\def\labelstyle{\displaystyle} {E_1}:\quad
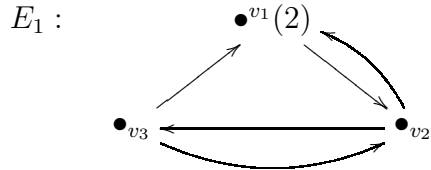
\xymatrix{ {} & \bullet^{v_1}  \ar[rd] {(2)}  &  {} \\
\bullet_{v_3}  \ar[ru] \ar@/^{-15pt}/ [rr]&  & \bullet_{v_2}
\ar[ll] \ar@/^{-10pt}/ [lu] \\
}}$$ \vspace{.2truecm}

(here $(2)$ means that there are two edges from $v_1$ to $v_2$). Then, $E_1=E_2(v_1\hookrightarrow
v_3)$, where
$${
\def\labelstyle{\displaystyle} {E_2}:\quad
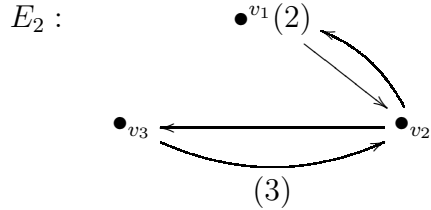
\xymatrix{ {} & \bullet^{v_1}  \ar[rd] {(2)}  &  {} \\
\bullet_{v_3}  \ar@/^{-15pt}/ [rr]_{(3)}&  & \bullet_{v_2}
\ar[ll] \ar@/^{-10pt}/ [lu] \\
}}$$ \vspace{.2truecm}

Also, $E_2=[E_3(v_3\hookrightarrow v_2)](v_1\hookrightarrow v_2)$, where
$${
\def\labelstyle{\displaystyle} {E_3}:\quad
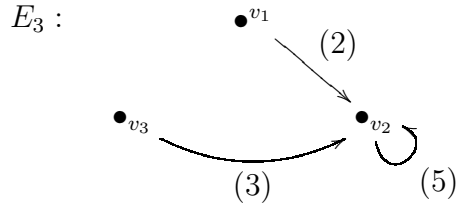
\xymatrix{ {} & \bullet^{v_1}  \ar[rd]^{(2)}  &  {} \\
\bullet_{v_3}  \ar@/^{-15pt}/ [rr]_{(3)}&  & \bullet_{v_2} \ar@(d,r)_{(5)}
\\
}}$$ \vspace{.2truecm}

Partitioning the set of edges emitted by $v_1$ and $v_2$ in singletons, we get
$${
\def\labelstyle{\displaystyle} {E_4}:\quad
\xymatrix{{\bullet}^{v_1} \ar[drr] &{\bullet}^{v_2} \ar[dr] & {\bullet}^{v_3} \ar[d] &
{\bullet}^{v_4} \ar[dl] & {\bullet}^{v_5} \ar[dll]\\ & & {\bullet}_{w} \ar@(dl,dr)_{(5)}  & & }} $$
\vspace{.2truecm}

Thus, $E_5=[[[[E_4(v_1\hookrightarrow w)](v_2\hookrightarrow w)](v_3\hookrightarrow
w)](v_4\hookrightarrow w)](v_5\hookrightarrow w)$ is the graph
$${
\def\labelstyle{\displaystyle} {E_5}:\quad
\xymatrix{ & & \bullet^{v_1} \ar@/^{-10pt}/ [dd]& & \\
 \bullet^{v_5} \ar@/^{-10pt}/ [rrd] & &  & & \bullet^{v_2} \ar@/^{-10pt}/ [lld] \\
  & & \bullet^{w} \ar@/^{-10pt}/ [uu]  \ar@/^{-10pt}/ [llu] \ar@/^{-10pt}/ [rru] \ar@/^{-10pt}/ [rdd]\ar@/^{-10pt}/ [ldd]& & \\
   & &  & & \\
    &\bullet^{v_4}\ar@/^{-10pt}/ [ruu] &  & \bullet^{v_3}\ar@/^{-10pt}/ [luu]& \\
}} $$ \vspace{.2truecm}

Recall that ${\rm M}_2(L_5)\cong L(B_5^2)$, where
$${
\def\labelstyle{\displaystyle} {B_5^2}:\quad
\xymatrix{ \bullet^{v_1}\ar@(ul,dl)_{(5)} & \bullet^{v_2} \ar[l]}}.
$$ \vspace{.2truecm}

Consider $\widehat{B_5^2}$
$${
\def\labelstyle{\displaystyle} \widehat{B_5^2}:\quad
\xymatrix{ & & \bullet^{v_1}\uloopr{} \ar[drr] \ar[dddr] \ar[dddl] \ar@/^{-10pt}/[dll]& & \\
\bullet^{v_5} \ar@(ul,dl)\ar[urr] \ar@/^{-10pt}/[rrrr]\ar@/^{-10pt}/[ddr]\ar@/^{-10pt}/[ddrrr]& &  & & \bullet^{v_2}\ar@(ur,dr) \ar[ddl] \ar[ddlll] \ar[llll] \ar@/^{-10pt}/[ull] \\
   & & \bullet^{w}\ar[uu]\ar[urr]\ar[ull]\ar[dr]\ar[dl] & & \\
    &\bullet^{v_4}\ar@(d,l)\ar[uul] \ar@/^{-10pt}/[rr]\ar@/^{-10pt}/[uul]\ar@/^{-10pt}/[uuur] &  & \bullet^{v_3}\ar@(d,r) \ar[ll] \ar@/^{-10pt}/[uur] \ar[uulll] \ar@/^{-10pt}/[uuul]& }} $$ \vspace{.2truecm}

We have $L(B_5^2)\cong L(\widehat{B_5^2})$ by Corollary \ref{isomaxoutsplit}. Then, we get
$F_1=[[[[[\widehat{B_5^2}](v_1\hookrightarrow v_2)](v_1\hookrightarrow v_3)](v_1\hookrightarrow
v_4)](v_1\hookrightarrow v_5)](v_1\hookrightarrow w)$
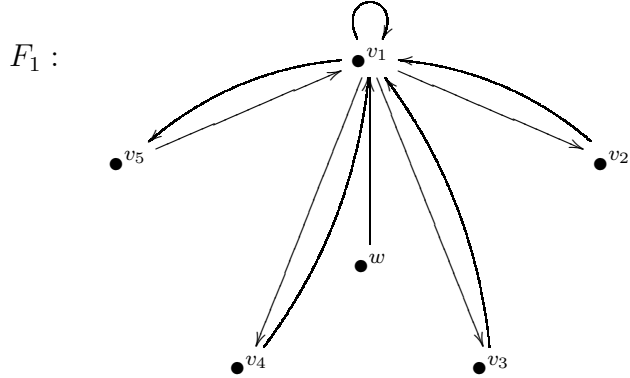
$${
\def\labelstyle{\displaystyle} {F_1}:\quad
\xymatrix{ & & \bullet^{v_1}\uloopr{} \ar[drr] \ar[dddr] \ar[dddl] \ar@/^{-10pt}/[dll]& & \\
\bullet^{v_5} \ar[urr] & &  & & \bullet^{v_2} \ar@/^{-10pt}/[ull] \\
   & & \bullet^{w}\ar[uu] & & \\
    &\bullet^{v_4}\ar@/^{-10pt}/[uuur] &  & \bullet^{v_3} \ar@/^{-10pt}/[uuul]& }}$$
and hence $L(\widehat{B_5^2})\cong L(F_1)$ by several applications of Theorem \ref{meu1_gen}.
Finally, $E_5=F_1(w\hookrightarrow v_1)$. Thus, $L(E_5)\cong L(F_1)$ by Theorem \ref{meu1_gen}, as
desired.}
\end{point}

\bigskip

In the final two subsections we analyze the remaining three graphs.  Since the $K_0$-data of these
graphs is not of the form $(\Z/(n-1)\Z,\overline{k})$, connections between the Leavitt path
algebras of these three graphs and algebras of the form ${\rm M}_k(L_n)$ are not of issue.

\begin{point}{\rm \textbf{$(K_0(L(E)), [1_{L(E)}])\cong (\mathbb{Z}/2\mathbb{Z}\oplus
\mathbb{Z}/2\mathbb{Z}, (\overline{0}, \overline{0}))$}: In this situation we have $1$ graph,
listed as follows:
$${E_{1}^6}: \quad {
\def\labelstyle{\displaystyle}
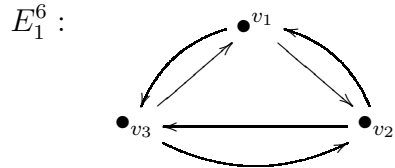
\xymatrix{ {} & \bullet^{v_1}  \ar[rd] \ar@/^{-10pt}/ [ld] &  {} \\
\bullet_{v_3}  \ar[ru] \ar@/^{-15pt}/ [rr]&  & \bullet_{v_2}
 \ar[ll]
\ar@/^{-10pt}/ [lu] \\
}}
$$
\vspace{.2truecm}}
\end{point}

\bigskip

\begin{point}{\rm \textbf{$(K_0(L(E)), [1_{L(E)}])\cong (\mathbb{Z}, 0)$}: In this
situation we have $2$ graphs, listed as follows:
$${E_{1}^7}: \quad {
\def\labelstyle{\displaystyle}
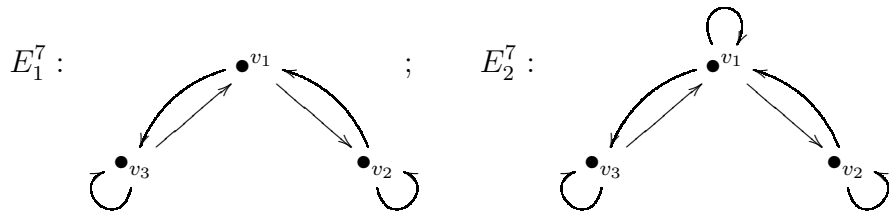
\xymatrix{ {} & \bullet^{v_1}  \ar[rd] \ar@/^{-10pt}/ [ld] &  {} \\
\bullet_{v_3} \ar@(d,l) \ar[ru] &  & \bullet_{v_2} \ar@(d,r)
\ar@/^{-10pt}/ [lu] \\
}}; \qquad {E_{2}^7}: \quad {
\def\labelstyle{\displaystyle}
\xymatrix{ {} & \bullet^{v_1} \uloopr{} \ar[rd] \ar@/^{-10pt}/ [ld] &  {} \\
\bullet_{v_3} \ar@(d,l) \ar[ru] &  & \bullet_{v_2} \ar@(d,r)
\ar@/^{-10pt}/ [lu] \\
}}
$$
\vspace{.2truecm}

By noticing that $E_{2}^7(v_3\hookrightarrow v_1)=E_{1}^7$, these Leavitt path algebras are
isomorphic by Theorem \ref{meu1_gen}, so we are done.}
\end{point}

\medskip

Thus we have answered in the affirmative another specific case of The Classification Question for
purely infinite simple unital Leavitt path algebras.

\begin{prop}\label{threevertices}
Suppose $E$ and $F$ are graphs satisfying Condition {\rm (Sing)} for which $L(E)$ and $L(F)$ are
purely infinite simple unital, and $|E^0| = |F^0| = 3$. If $K_0(L(E))\cong K_0(L(F))$ via an
isomorphism $\phi$ for which $\phi([1_{L(E)}])=[1_{L(F)}]$, then $L(E)\cong L(F)$.  Moreover, for
any such graph $E$ for which $(K_0(L(E)), [1_{L(E)}])\cong (\Z/(n-1)\Z,\overline{k})$, then in fact
we have $L(E)\cong {\rm M}_k(L_n)$.
\end{prop}

\section{Isomorphisms between Leavitt algebras and their matrices}

In this final section we deal with the problem of determining values of $k$ and $n$ for which
$L_n\cong {\rm M}_k(L_n)$. This problem was completely solved by three of the authors in
\cite{AbAnhP} using combinatorial arguments, where we show that the necessary and sufficient
condition for such an isomorphism is that $\mbox{g.c.d.}(k, n-1)=1$.  As it turns out, one
direction of this implication was already established in \cite{L}, where Leavitt shows that $L_n$
and ${\rm M}_k(L_n)$ are not isomorphic whenever $\mbox{g.c.d.}(k, n-1)>1$.  Thus throughout this
section we consider only situations in which $\mbox{g.c.d.}(k, n-1)=1$. Our goal here is to prove
the isomorphism $L_n\cong {\rm M}_k(L_n)$, using arguments afforded by  Theorems \ref{meu1_gen} and
\ref{outsplitiso}. The graph and combinatorial approaches to this isomorphism question are
essentially independent.

First, fix graphs
$${
\def\labelstyle{\displaystyle} {R_n^k}:\quad
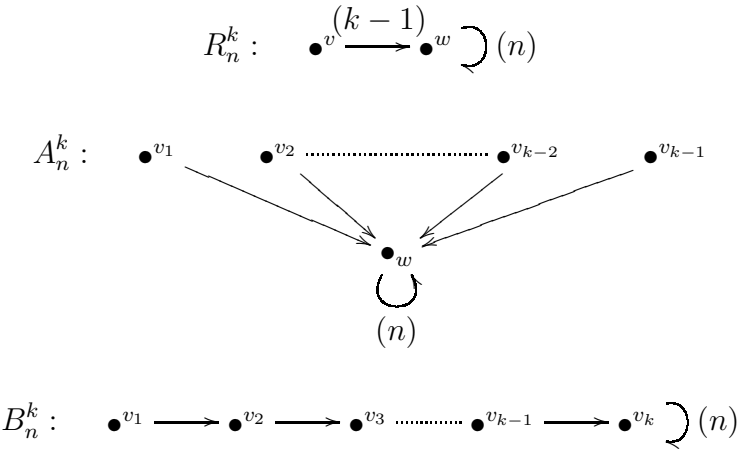
\xymatrix{{\bullet}^{v} \ar [r] ^{(k-1)} & {\bullet}^{w} \ar@(ur,dr)^{(n)}} }$$

$${
\def\labelstyle{\displaystyle} {A_n^k}:\quad
\xymatrix{{\bullet}^{v_1}  \ar[drr] & {\bullet}^{v_2}  \ar[dr] \ar@{.}[rr] &  & {\bullet}^{v_{k-2}}
\ar[dl]&{\bullet}^{v_{k-1}} \ar[dll] \\& & {\bullet}_w \ar@(dl,dr)_{(n)}
 & & }}$$

$${
\def\labelstyle{\displaystyle} {B_n^k}:\quad
\xymatrix{{\bullet}^{v_1} \ar [r]  & {\bullet}^{v_2} \ar [r] &
 {\bullet}^{v_3} \ar@{.}[r] &
{\bullet}^{v_{k-1}} \ar [r]  & {\bullet}^{v_k} \ar@(ur,dr)^{(n)}} }$$

Then, we have the following result.

\begin{lem}\label{L:OpenTails}
For each $n\geq 1$, the algebras $L(A_n^k)$, $L(B_n^k)$ and $L(R_n^k)$ are isomorphic. Moreover,
each of these algebras
 is isomorphic to the matrix algebra ${\rm M}_k(L_n)$.
\end{lem}
\begin{proof}
It is clear that $A_n^k$ is obtained from $R_n^k$ by splitting $s^{-1}(v)$ in $k-1$ singletons,
whence $L(R_n^k)\cong L(A_n^k)$ by Theorem \ref{outsplitiso}. On the other side, rewrite $A_n^k$ as
$${
\def\labelstyle{\displaystyle} {B_{k-1}^1}:\quad
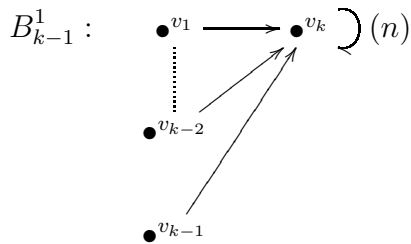
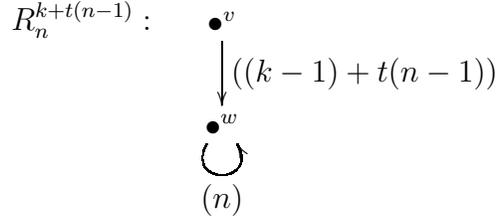
\xymatrix{{\bullet}^{v_1} \ar [r]  \ar@{.}[d] & {\bullet}^{v_k}
\ar@(ur,dr)^{(n)} \\
{\bullet}^{v_{k-2}}\ar [ur]  & \\
{\bullet}^{v_{k-1}}\ar [uur]  & \\}}$$ Then, $B_{k-2}^2=B_{k-1}^1(v_2\hookrightarrow v_1)$, where
$${
\def\labelstyle{\displaystyle} {B_{k-2}^2}:\quad
\xymatrix{{\bullet}^{v_1} \ar [r]   & {\bullet}^{v_2} \ar [r]
\ar@{.}[d] & {\bullet}^{v_k} \ar@(ur,dr)^{(n)}\\
 & {\bullet}^{v_{k-2}}\ar [ur]  & \\
 & {\bullet}^{v_{k-1}}\ar [uur]  & \\}}$$
Hence, $L(B_{k-1}^1)\cong L(B_{k-2}^2)$ by Theorem \ref{meu1_gen}. (Note that Condition (L) holds
trivially.) Recurrence on this argument produces a chain of graphs whose Leavitt path algebras are
isomorphic. This chain ends in
$${
\def\labelstyle{\displaystyle} {B_2^{k-2}}:\quad
\xymatrix{{\bullet}^{v_1} \ar [r] & {\bullet}^{v_2} \ar [r]  & {\bullet}^{v_3} \ar@{.}[r] &
{\bullet}^{v_{k-3}} \ar [r]  & {\bullet}^{v_{k}}
\ar@(ur,dr)^{(n)}\\
 & & & {\bullet}^{v_{k-2}}\ar [ur]  & \\
  & & & {\bullet}^{v_{k-1}}\ar [uur]  & \\}}$$

$${
\def\labelstyle{\displaystyle} {B_1^{k-1}}:\quad
\xymatrix{{\bullet}^{v_1} \ar [r]  & {\bullet}^{v_2} \ar [r] & {\bullet}^{v_3} \ar@{.}[r] &
{\bullet}^{v_{k-2}} \ar [r]  & {\bullet}^{v_{k}}
\ar@(ur,dr)^{(n)}\\
 & & & {\bullet}^{v_{k-1}}\ar [ur]  & \\}}$$

with the final step being
$${
\def\labelstyle{\displaystyle} {B_n^k}:\quad
\xymatrix{{\bullet}^{v_1} \ar [r]  & {\bullet}^{v_2} \ar [r] &
 {\bullet}^{v_3} \ar@{.}[r] &
{\bullet}^{v_{k-1}} \ar [r]  & {\bullet}^{v_k} \ar@(ur,dr)^{(n)}} }.$$
 This yields the asserted isomorphisms between the three indicated Leavitt path algebras.  The
final statement
 follows immediately from \cite[Proposition 13]{AA2}, in which the isomorphism $L(B_n^k) \cong {\rm M}_k(L_n)$ is established.
\end{proof}

\medskip

Since any graph having Leavitt path algebra isomorphic to ${\rm M}_k(L_n)$ must satisfy Condition
(L) (as any such algebra is purely infinite simple), Theorem \ref{meu1_gen} may be invoked in all
situations throughout the sequel. Specifically, all of the isomorphism results of this section
 hold regardless of the cardinality of the field $K$.

Here is our first result about isomorphisms between matrix rings over Leavitt algebras.

\begin{prop}\label{L:congruent mod n-1}
For every $t\geq 0$, for every $k\geq 1$ and for every $n\geq 2$,
$${\rm M}_k(L_n)\cong {\rm M}_{k+t(n-1)}(L_n).$$
\end{prop}
\begin{proof}
We will prove the result by induction on $t$. The case $t=0$ being clear, we suppose than the
result holds for $t-1$. By \cite[Proposition 13]{AA2} we have  ${\rm M}_{k+t(n-1)}(L_n)\cong
L(R_n^{k+t(n-1)})$, where

$${
\def\labelstyle{\displaystyle} {R_n^{k+t(n-1)}}:\quad
\xymatrix{{\bullet}^{v} \ar[d] ^{((k-1)+t(n-1))}\\  {\bullet}^w \ar@(dl,dr)_{(n)} }}
$$ \vspace{.2truecm}

Splitting the edges emitted by $v$ in two sets, one with $(k-1)+t(n-1)$ edges, and the other with
$(n-1)$ edges, we get $E_1=(R_n^{k+t(n-1)})_s(\mathcal{P})$

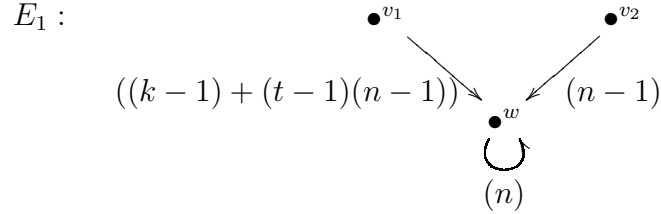
$${
\def\labelstyle{\displaystyle} {E_1}:\quad
\xymatrix{{\bullet}^{v_1} \ar[dr]_{((k-1)+(t-1)(n-1))} & & {\bullet}^{v_2} \ar[dl] ^{(n-1)}\\
&{\bullet}^w \ar@(dl,dr)_{(n)} & }}
$$ \vspace{.2truecm}

and $L(R_n^{k+t(n-1)})\cong L(E_1)$ by Theorem \ref{outsplitiso}. Now, Let
$E_2=E_1(v_2\hookrightarrow w)$,
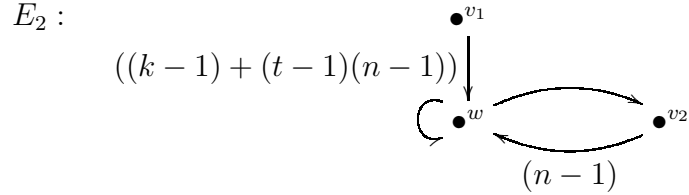
$${
\def\labelstyle{\displaystyle}E_2:\quad \xymatrix{\bullet^{v_1}
\ar[d]_{((k-1)+(t-1)(n-1))} & & \\ {\bullet}^{w} \ar@(ul,dl) \ar@/^1pc/ [rr] & & {\bullet}^{v_2}
\ar@/^1pc/ [ll] ^{(n-1)} }}$$ \vspace{.2truecm}

By Theorem \ref{meu1_gen}, $L(E_1)\cong L(E_2)$. Take $E_3$ to be
$${
\def\labelstyle{\displaystyle}E_3:\quad 
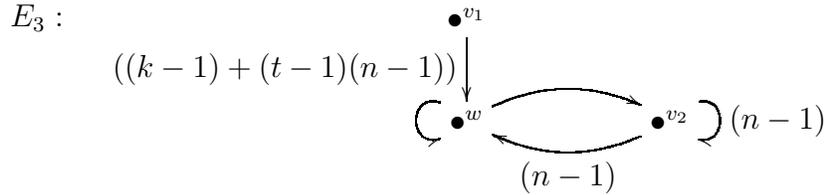
\xymatrix{\bullet^{v_1}
\ar[d]_{((k-1)+(t-1)(n-1))} & &  \\ {\bullet}^{w} \ar@(ul,dl) \ar@/^1pc/ [rr] & & {\bullet}^{v_2}
\ar@(ur,dr)^{(n-1)} \ar@/^1pc/ [ll] ^{(n-1)} }}$$ \vspace{.2truecm} and notice that $E_2=[\cdots
[E_3(w\hookrightarrow v_2)](w\hookrightarrow v_2)]\cdots ](w\hookrightarrow v_2)$, $(n-1)$ times,
so that $L(E_2)\cong L(E_3)$ by Theorem \ref{meu1_gen}. Now, take $E_4$
$${
\def\labelstyle{\displaystyle}E_4:\quad 
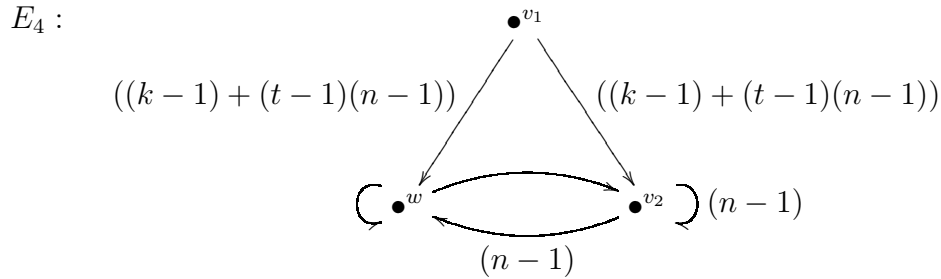
\xymatrix{ & \bullet^{v_1}
\ar[ddl]_{((k-1)+(t-1)(n-1))} \ar[ddr]^{((k-1)+(t-1)(n-1))}& \\ & & \\
{\bullet}^{w} \ar@(ul,dl) \ar@/^1pc/ [rr] & & {\bullet}^{v_2} \ar@(ur,dr)^{(n-1)} \ar@/^1pc/ [ll]
^{(n-1)} }}$$ \vspace{.2truecm} and notice that $E_3=[\cdots [E_4(w\hookrightarrow
v_1)](w\hookrightarrow v_1)]\cdots](w\hookrightarrow v_1)$, $(k-1)+(t-1)(n-1)$ times, so that
$L(E_3)\cong L(E_4)$ by Theorem \ref{meu1_gen}. Finally, take $R_n^{k+(t-1)(n-1)}$
$${
\def\labelstyle{\displaystyle} {R_n^{k+(t-1)(n-1)}}:\quad
\xymatrix{{\bullet}^{v} \ar[d] ^{((k-1)+(t-1)(n-1))}\\  {\bullet}^w \ar@(dl,dr)_{(n)} }}
$$ \vspace{.2truecm}
and notice that $E_4$ is the out-splitting of $R_n^{k+(t-1)(n-1)}$ over the edges emitted by $w$ in
two sets, one with $(n-1)$ edges, and the other a singleton. Thus, $L(E_4)\cong
L(R_n^{k+(t-1)(n-1)})$; but the latter is isomorphic to ${\rm M}_k(L_n)$ by the induction
hypothesis.  Thus ${\rm M}_{k+t(n-1)}(L_n)\cong {\rm M}_k(L_n)$, which completes the induction
step.
\end{proof}

We note that the conclusion of Proposition \ref{L:congruent mod n-1} also follows from the fact
that the free left $L_n$-modules of ranks $k$ and $k+t(n-1)$ are isomorphic, so that the
endomorphism rings of these modules are isomorphic, and such endomorphism rings are in turn
isomorphic to the indicated matrix rings.

Here is our second result about isomorphisms between matrix rings over Leavitt algebras.

\begin{prop}\label{P:dividingpower}
Let $n\geq 2, k\geq 2$ be such that $k$ divides $n$. Then $L_n\cong {\rm M}_k(L_n)$.
\end{prop}
\begin{proof}
We have that $n=kl$ for some $l$.  We recall again that ${\rm M}_k(L_n)\cong L(R_n^{k})$,
$${
\def\labelstyle{\displaystyle} {R_n^{k}}:\quad
\xymatrix{{\bullet}^{v} \ar[d] ^{(k-1)}\\
{\bullet}^w \ar@(dl,dr)_{(n)} }}
$$ \vspace{.2truecm}

Of course $n=(k-1)l+l$. Consider $E_1=[\cdots [R_n^{k}(v\hookrightarrow w)](v\hookrightarrow
w)]\cdots ](v\hookrightarrow w)$, $l$ times
$${
\def\labelstyle{\displaystyle}E_1:\quad \xymatrix{  {\bullet}^{w} \ar@(ul,dl)_{(l)} \ar@/^1pc/
[rr]^{(l)} & & {\bullet}^{v}  \ar@/^1pc/ [ll] ^{(k-1)} }}$$ \vspace{.2truecm}

By Theorem \ref{meu1_gen}, $L(R_n^{k})\cong L(E_1)$. Now, take $E_2$
$${
\def\labelstyle{\displaystyle}E_2:\quad \xymatrix{  {\bullet}^{w} \ar@(ul,dl)_{(l)} \ar@/^1pc/
[rr]^{(l)} & & {\bullet}^{v} \ar@(ur,dr)^{(k-1)l} \ar@/^1pc/ [ll] ^{(k-1)l} }}$$ \vspace{.2truecm}
and notice that $E_1=[\cdots [E_2(w\hookrightarrow v)](w\hookrightarrow v)]\cdots
](w\hookrightarrow v)$, $k-1$ times, so that $L(E_1)\cong L(E_2)$ by Theorem \ref{meu1_gen}.
Finally, consider $R_n$
$${
\def\labelstyle{\displaystyle}R_n:\quad
\xymatrix{{\bullet^w}\ar@(ur,dr)^{(n)} }}$$\vspace{.2truecm}

If we consider a partition of the edges emitted by $w$ in two sets with $l$ and $(k-1)l$ edges
respectively, then $(R_n)_s(\mathcal{P})=E_2$, so that $L(R_n)\cong L(E_2)$ by Theorem
\ref{outsplitiso}. As $L_n\cong L(R_n)$, the desired result holds.
\end{proof}

\bigskip

The final goal of this article is to use our two ``Change the Graph" isomorphisms to establish the
isomorphism $L_n\cong {\rm M}_k(L_n)$ whenever ${\rm g.c.d.}(k,n-1)=1$. As mentioned
 previously, this isomorphism was established in \cite{AbAnhP}, using completely different techniques.  As a consequence of the current discussion,
 we obtain yet more evidence suggesting an affirmative answer to The Classification Question.

 As we shall see, establishing the isomorphism $L_n\cong {\rm M}_k(L_n)$ utilizes a ten step process, where each step requires the use of one or the
  other of the two Change the Graph theorems.  Steps 1 through 3, and 5 through 10, are relatively transparent; however,  Step 4 requires some additional work, which we
  take care of in the next few results.

  We begin by relating Theorem \ref{meu1_gen} to a matrix operator.

  \begin{defi}\label{Phi}   {\rm Let $M\in {\rm M}_p(\Z)$.  Define $\Phi_M: {\rm M}_p(\Z) \rightarrow {\rm M}_p(\Z)$ by setting, for each $A\in {\rm M}_p(\Z)$,
$$ \Phi_M(A)=MA+(I_p-M)$$
where $ I_p$ is the identity matrix in ${\rm M}_p(\Z^+)$.  It is easy to check that for $M_1, M_2,
..., M_t \in {\rm M}_p(\Z)$ we have $\Phi_{M_t}\circ \cdots \circ \Phi_{M_2}\circ
\Phi_{M_1}=\Phi_{M_t \cdots M_2M_1}$.}
\end{defi}

Straightforward matrix arithmetic yields

\begin{lem}\label{PhiofAnonnegative}  Let $k,p\in \N$.  For
integers $1\leq s,t\leq p$ with $s\neq t$ let $K \in {\rm M}_p(\Z^+)$ denote the matrix
  $$K = I_p+ke_{st}$$
   where $e_{st}$ denotes the standard $(s,t)$
matrix unit in ${\rm M}_p(\Z^+)$. Suppose $A=(a_{ij})\in {\rm M}_p(\Z^+)$.

  \medskip

        (1)   If $a_{tt} \geq 1$ and $a_{st} \geq 1$ then
$\Phi_K(A)\in {\rm M}_p(\Z^+),$  and for each $1\leq j \leq p$ we have
        $(\Phi_K(A))_{sj}\geq (\Phi_K(A))_{tj}$ and
$(\Phi_K(A))_{tt}\geq
        1.$

        \medskip

        (2)   In addition,  if $a_{ts}\geq 1$ then
$(\Phi_K(A))_{ss}\geq 1$ and $(\Phi_K(A))_{ts}\geq
        1.$

  \end{lem}

The proof of the next result follows directly from Lemma \ref{PhiofAnonnegative} and the
construction presented in Theorem \ref{meu1_gen}.
\begin{corol}\label{LEisoLPhiE}
Let $k,p\in \N$.  Let $K \in {\rm M}_p(\Z^+)$ denote the matrix $I_p+ke_{st}$ for some pair $1\leq
s,t\leq p$ with $s\neq t$. Suppose $A=(a_{ij})\in {\rm M}_p(\Z^+)$ has $a_{tt} \geq 1$ and $a_{st}
\geq 1$, and suppose that the associated graph $E_A$ satisfies Condition (L). Then
$$L(E_A)\cong L(E_{\Phi_K(A)}).$$
\end{corol}

  \begin{defi}\label{Smatrices}  {\rm  We identify some quantities which will be useful in the sequel.   Let $a,b$ be positive integers having
  ${\rm g.c.d.}(a,b)=1$.  Assume $a>b>1$.   We apply the standard Euclidean algorithm to find sequences of positive integers $r_0, r_1, ... ,r_m$ and $k_1, k_2, ... ,k_m$ for
   which
$$r_0=k_1r_1+r_2, \ \ \  r_1=k_2r_2+r_3,  \ \hdots , \ \ \ r_{m-2}=k_{m-1}r_{m-1}+r_m $$
where $r_0=a$, $r_1=b$, and $r_m={\rm g.c.d.}(a,b)=1$. (Note that $m\geq 2$ since $b>1$.)  It will
be notationally useful to add a non-standard additional equation to the end of this list by
defining
 $$r_{m-1}=k_m\cdot 1+r_{m+1}$$
(in other words, we set $k_m=r_{m-1}-1$ and $r_{m+1}=1$.)  We now define a collection of $3 \times
3$ matrices based on these sequences, by setting
$$S_0=\left(%
\begin{array}{ccc}
  1 & 0 & 0 \\
  0 & 1 & 0 \\
  0 & 0 & 1
\end{array}%
\right), \ \ \
S_1=\left(%
\begin{array}{ccc}
  1 & 0 & 0 \\
  k_1 & 1 & 0 \\
  0 & 0 & 1
\end{array}%
\right), $$
 and, for every integer $p$ with $1\leq p \leq m/2$,
$$S_{2p}=\left(%
\begin{array}{ccc}
  1 & k_{2p} & 0 \\
  0 & 1 & 0 \\
  0 & 0 & 1
\end{array}%
\right)\cdot S_{2p-1}, \ \ \
S_{2p+1}=\left(%
\begin{array}{ccc}
  1 & 0 & 0 \\
  k_{2p+1} & 1 & 0 \\
  0 & 0 & 1
\end{array}%
\right)\cdot S_{2p}.
$$
We note that $S_i\in {\rm M}_p(\Z^+)$ and ${\rm det}(S_i)=1$ for all $0\leq i \leq m$. }
\end{defi}

\begin{prop}\label{L:Graph Matrix}  Let $a>b$ be positive integers having
  ${\rm g.c.d.}(a,b)=1$.  Suppose also $b>1$.   We use the notation as in Definition \ref{Smatrices}.  We denote the specific matrix $S_m$ by
$$S_{m}=\left(%
\begin{array}{ccc}
  x_1 & y_1 & 0 \\
  x_2 & y_2 & 0 \\
  0 & 0 & 1
\end{array}%
\right).$$ In particular we have $y_2\geq 1$. Then

(1)  $x_1b-y_1a=1$ and $x_2b-y_2a=-1$.

 (2)  $x_1 - y_1 \geq 1$ and $x_2 - y_2 \geq 0$.

 (3)  $x_1 + x_2 = a$ and $y_1 + y_2 = b$.

 \end{prop}

\begin{proof}

Using the matrix equations given in Definition \ref{Smatrices}, an easy induction argument shows,
for each integer $p$ having $1\leq p \leq m/2$, that
$$ S_{2p}\left(%
\begin{array}{c}
  r_1 \\
  -r_0 \\
  0\\
\end{array}%
\right)=\left(%
\begin{array}{c}
  r_{2p+1} \\
  -r_{2p}\\
  0\\
\end{array}%
\right) \ \
 \mbox{ and } \ \
S_{2p+1}\left(%
\begin{array}{c}
  r_1 \\
  -r_0 \\
  0 \\
\end{array}%
\right)=\left(%
\begin{array}{c}
  r_{2p+1} \\
  -r_{2p+2}\\
  0 \\
\end{array}%
\right).  $$
 Thus in particular, using that $r_m=r_{m+1}=1$, we get
$$S_m\left(%
\begin{array}{c}
  r_1 \\
  -r_0 \\
  0\\
\end{array}%
\right)=\left(%
\begin{array}{c}
 1 \\
 -1\\
 0\\
\end{array}%
\right).$$
 But recall that by definition we have $r_1=b$ and $r_0=a$, so the previous equation becomes
$$\left(%
\begin{array}{ccc}
  x_1 & y_1 & 0 \\
  x_2 & y_2 & 0 \\
  0 & 0 & 1
\end{array}%
\right)
\left(%
\begin{array}{c}
  b \\
  -a \\
  0\\
\end{array}%
\right)=\left(%
\begin{array}{c}
 1 \\
 -1\\
 0\\
\end{array}%
\right).$$
 This matrix equation yields (1).   The equation $x_1 - y_1 \geq 1$ of (2) follows immediately from (1) and the
  hypothesis that $a>b$.   For the other part of (2), note that $-1 = x_2b - y_2a < x_2b-y_2b$ (since $a>b>0$ and $y_2 > 0$), so $-1< (x_2-y_2)b \in \Z$, so
  necessarily $x_2-y_2\geq 0$.

   Now recall that ${\rm det}(S_i)=1$ for all $0\leq i \leq m$, so in particular we have
 ${\rm det}(S_m)=1$, so that $x_1y_2-y_1x_2=1$.  We incorporate the previous two pieces of information in the single matrix equation
$$\left(%
\begin{array}{ccc}
  x_1 & y_1 & 0 \\
  x_2 & y_2 & 0 \\
  0 & 0 & 1
\end{array}%
\right)
\left(%
\begin{array}{ccc}
  b & -y_1 & 0 \\
  -a & x_1 & 0 \\
  0 & 0 & 1
\end{array}%
\right) =
\left(%
\begin{array}{ccc}
  1 & 0 & 0 \\
  -1& 1 & 0 \\
  0 & 0 & 1
\end{array}%
\right).$$
 But then by using this same information it is easy to check that

$$
\left(%
\begin{array}{ccc}
  b & -y_1 & 0 \\
  -a & x_1 & 0 \\
  0 & 0 & 1
\end{array}%
\right)^{-1} =
\left(%
\begin{array}{ccc}
  x_1 & y_1 & 0 \\
  a & b & 0 \\
  0 & 0 & 1
\end{array}%
\right),$$
 so that right multiplying gives

$$\left(%
\begin{array}{ccc}
  x_1 & y_1 & 0 \\
  x_2 & y_2 & 0 \\
  0 & 0 & 1
\end{array}%
\right) =
\left(%
\begin{array}{ccc}
  1 & 0 & 0 \\
  -1& 1 & 0 \\
  0 & 0 & 1
\end{array}%
\right)
\left(%
\begin{array}{ccc}
  x_1 & y_1 & 0 \\
  a & b & 0 \\
  0 & 0 & 1
\end{array}%
\right) =
\left(%
\begin{array}{ccc}
  x_1 & y_1 & 0 \\
  a-x_1 & b-y_1 & 0 \\
  0 & 0 & 1
\end{array}%
\right).$$
 In particular $x_2 = a-x_1$ and $y_2=b-y_1$, from which (3) follows immediately.

\end{proof}

Were our aim solely to achieve a number-theoretic result regarding the Euclidean algorithm, we
would have constructed analogous $S$-matrices inside ${\rm M}_2(\Z^+)$ rather than inside ${\rm
M}_3(\Z^+)$,
 since clearly all of the germane computations take place in the upper $2\times 2$ blocks of the $S_i$.  However, as we shall see below, our one application
 of the combinatorial facts provided in Proposition \ref{L:Graph Matrix} will be in the context of $3\times 3$ matrices, so we choose to formulate the result
 accordingly.

We are now in position to demonstrate our main isomorphism result, a result which will lead to
verification of another piece of The Classification Question.

\begin{theor}\label{Th:Medium Fish}
Let $n,d$ be positive integers having $\mbox{g.c.d.}(d,n-1)=1$, and let $K$ be any field.   Then
$L_n\cong\emph{M}_d(L_n).$
\end{theor}
\begin{proof}
The result is trivial for $d=1$. Now suppose $d=2$.   Then either $n$ is even, whence $d$ divides
$n$ and so the result holds by Proposition  \ref{P:dividingpower}, or $n$ is odd, whence
$\mbox{g.c.d.}(d,n-1)=2$, contradicting the hypothesis.  So we may assume that $d\geq 3$.  Also, by
Proposition \ref{L:congruent mod n-1} we can assume that $d\leq n-2$.  Now write $n=dt+r$ with
$0\leq r \leq d-1.$  (In particular,  $d-r+1 \geq 2$.) If $r=0$ then $d$ divides $n$, so that the
result holds by Proposition \ref{P:dividingpower}. If $r=1$ then $d=\mbox{g.c.d.}(d, n-1)$,
contradicting the hypothesis.  So we may also assume that $r\geq 2$.  In particular we have
$t+r-1\geq 1$ for all $t\geq 0$.

It will be clear that each of the graphs encountered in this proof satisfies Condition (L).  Thus
Theorem \ref{meu1_gen} may be invoked throughout, without regard to the size of the field of
scalars.

\medskip

  Recall that $L_n \cong L(R_n)$ where $R_n$ is the graph
$$R_n: \ \ \ \ \ \xymatrix{
  \bullet^{v} \ar@(ul,dl)[]_{(n)}},$$
while ${\rm M}_d(L_n)\cong L(R_n^d)$ where $R_n^d$ is the graph
$$ R_n^d: \ \ \ \
\ \xymatrix{
    \bullet^{v_1} \ar@(ul,dl)[]_{(n)} &
   \bullet^{v_2} \ar[l]_{(d-1)}
  }$$
We establish the desired result by building a ten step sequence of isomorphisms which starts with
$L(R_n^d)$ and ends with $L(R_n)$.

\medskip

  {\bf Step 1.}  Consider
$E_1=R_n^d(\stackrel{t\mbox{ times}}{v_2\hookrightarrow
  v_1})$,
  $$E_1: \ \ \ \ \xymatrix{
    \bullet^{v_1} \ar@(ul,dl)[]_{(t+r)} \ar@/^/[r]^{(t)} &
     \bullet^{v_2}  \ar@/^/[l]^{(d-1)} }$$
Then
$$L(R_n^d)\cong L(E_1) \mbox{ by Theorem \ref{meu1_gen}.}$$

\medskip

 {\bf Step 2.}   Splitting the set of edges of $E_1$ emitted by $v_2$ in two sets of $1$ and $d-2$ edges respectively, we get $E_2=(E_1)_s(\emph{P})$
$$E_2: \xymatrix{
& & \bullet^{v_2}  \ar@/^/[dll]^{(1)} \\
\bullet^{v_1} \ar@(ul,dl)[]_{(t+r)} \ar@/^/[urr]^{(t)}
\ar@/_/[drr]_{(t)} & &\\
 & & \bullet^{v_3}  \ar@/_/[ull]_{(d-2)}
 }$$
whence
$$L(E_1)\cong L(E_2) \mbox{ by Theorem \ref{outsplitiso}.}$$
(Note that $d\geq 3$ guarantees that the quantity $d-2$ is nonnegative.)

\medskip

{\bf Step 3.}   Consider $E_3=E_2(v_2\hookrightarrow
  v_1)$
  $$E_3: \xymatrix{
& & \bullet^{v_2}  \ar@/^/[dll]^{(1)} \\
\bullet^{v_1} \ar@(ul,dl)[]_{(t+r-1)} \ar@/^/[urr]^{(t+1)}
\ar@/_/[drr]_{(t)} & & \\
 & & \bullet^{v_3}  \ar@/_/[ull]_{(d-2)}
 }$$
Then we get
$$L(E_2)\cong L(E_3) \mbox{ by Theorem \ref{meu1_gen}.}$$

\medskip

{\bf Step 4.}  We are now in position to use the number-theoretic results described above.    Let
$A$ denote the matrix $A_{E_3}$ of the graph $E_3$; that is,
$$A=\left(%
\begin{array}{ccc}
  t+r-1 & t+1 & t\\
  1 & 0 & 0\\
  d-2 & 0 & 0
\end{array}%
\right),$$
 so that $E_A=E_3$.

Since $\mbox{g.c.d.}(d,r-1)=1$ we have $\mbox{g.c.d.}(d,d-r+1)=1.$ Let $r_0=a=d$, $r_1=b=d-r+1$.
Note that $a>b$ (since $r\geq 2$)
 and $b>1$ (since $r\leq d-1$). So we may apply the analysis given in Definition \ref{Smatrices} to the pair $a=d,b=d-r+1$ to produce the
  indicated matrix $S_m$.

  Now define, for $1\leq i \leq m$,
  $$K_i = I_3 + k_ie_{21} \mbox{ for } i \mbox{ odd,} \ \ \ \mbox{ and } \ \ \  K_i = I_3 + k_ie_{12} \mbox{ for } i \mbox{ even.}$$
  Then by construction we have $S_m = K_m\cdots K_2K_1$.

  Because $K_1 = I_3 + k_1e_{21}$ and $a_{11}= t+r-1\geq 1$ and $a_{21}=1$, we have by Lemma \ref{PhiofAnonnegative}(1) that $\Phi_{K_1}(A) \in {\rm M}_3(\Z^+)$.
   Moreover, since $a_{12}=t+1\geq 1$, Lemma \ref{PhiofAnonnegative}(2) yields $(\Phi_{K_1}(A))_{22}\geq 1$ and $(\Phi_{K_1}(A))_{12}\geq 1$.  Thus Lemma \ref{PhiofAnonnegative}
    may be applied at each step, and  we thereby conclude that the matrix $B$ defined by
   $$B=\Phi_{K_m}\circ\Phi_{K_{m-1}}\circ \cdots \circ \Phi_{K_1}(A)=\Phi_{S}(A)$$
   is  in ${\rm M}_p(\Z^+)$.  Now define $E_4 = E_B$ for the matrix $B$.   So we have
$$L(E_3)\cong L(E_B)=L(E_4) \mbox{ by an application of Corollary \ref{LEisoLPhiE} } m \mbox{ times.}$$

Prior to moving on to Step 5, we actually compute the values of the entries of the matrix
$B=(b_{ij})$.   By definition we have
$$B = \Phi_{S_m}(A) = \left(%
\begin{array}{ccc}
  x_1 & y_1 & 0\\
  x_2 & y_2 & 0\\
  0 & 0 & 1
\end{array}%
\right)
\left(%
\begin{array}{ccc}
  t+r-1 & t+1 & t\\
  1 & 0 & 0\\
  d-2 & 0 & 0
\end{array}%
\right)
+\left(%
\begin{array}{ccc}
  1-x_1 & -y_1 & 0\\
  -x_2 & 1-y_2 & 0\\
  0 & 0 & 0
\end{array}%
\right).$$
 So upon doing the matrix arithmetic we get
$$B = \left(%
\begin{array}{ccc}
  x_1(t+r-1)+y_1+1-x_1 & x_1(t+1)-y_1 & x_1t\\
  x_2(t+r-1)+y_2-x_2 & x_2(t+1)+1-y_2 & x_2t\\
  d-2 & 0 & 0
\end{array}%
\right)$$
 A pictorial description of $E_4=E_B$ is then given by

$$E_4: 
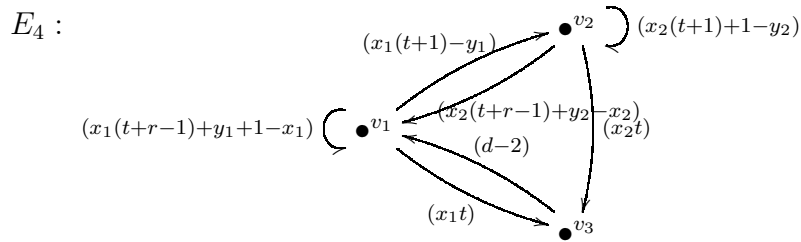
\xymatrix{
& & \bullet^{v_2} \ar@(ur,dr)[]^{(x_2(t+1)+1-y_2)} \ar@/^/[dll]^{(x_2(t+r-1)+y_2-x_2)} \ar@/^/[dd]^{(x_2t)} \\
\bullet^{v_1} \ar@(ul,dl)[]_{(x_1(t+r-1)+y_1+1-x_1)} \ar@/^/[urr]^{(x_1(t+1)-y_1)}
\ar@/_/[drr]_{(x_1t)} & & \\
 & & \bullet^{v_3}  \ar@/_/[ull]_{(d-2)}
 }$$


\medskip

{\bf Step 5.}    By Lemma \ref{L:Graph Matrix} we have $x_1-y_1 \geq 1$ and $x_2-y_2\geq 0$.  So we
may define the graph $E_5$ by setting
$$E_5=[E_4(\stackrel{(x_1-y_1)\mbox{ times}}{v_3\hookrightarrow
v_1})](\stackrel{(x_2-y_2)\mbox{ times}}{v_3\hookrightarrow v_2}).$$
 Using the information presented in the various parts of Lemma \ref{L:Graph Matrix},
it is tedious but straightforward to now verify each of the following equations.
$$\begin{array}{ll}
b_{11}-(x_1-y_1)b_{31}& =x_1(t+1)-y_1 \\
b_{12}-(x_1-y_1)b_{32}& =x_1(t+1)-y_1\\
b_{13}-(x_1-y_1)b_{33}+(x_1-y_1)&=x_1(t+1)-y_1\\
b_{21}-(x_2-y_2)b_{31}& =x_2(t+1)+1-y_2\\
b_{22}-(x_2-y_2)b_{32}& =x_2(t+1)+1-y_2\\
b_{23}-(x_2-y_2)b_{33}+(x_2-y_2) & =x_2(t+1)-y_2\\
\end{array}$$
For notational convenience we  define $n_1=x_1(t+1)-y_1$ and $n_2=x_2(t+1)+1-y_2$. Note that
$n_1\geq 1$ and $n_2\geq 1$ as a consequence of Lemma \ref{L:Graph Matrix}(2). Now using this list
of equations, we have that $E_5$ is the graph
 $$E_5: 
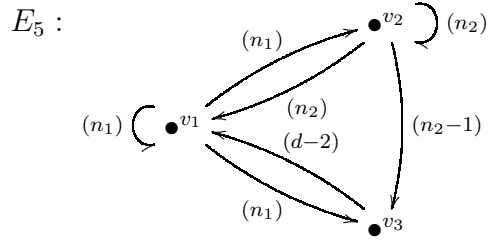
\xymatrix{
& & \bullet^{v_2} \ar@(ur,dr)[]^{(n_2)} \ar@/^/[dll]^{(n_2)} \ar@/^/[dd]^{(n_2-1)} \\
\bullet^{v_1} \ar@(ul,dl)[]_{(n_1)} \ar@/^/[urr]^{(n_1)}
\ar@/_/[drr]_{(n_1)} & & \\
 & & \bullet^{v_3}  \ar@/_/[ull]_{(d-2)}
 }$$
In particular,
$$L(E_4)\cong L(E_5) \mbox{ by Theorem \ref{meu1_gen}.}$$

\medskip

{\bf Step 6.}   Since
\begin{eqnarray*}
  n_1+n_2 &=& (x_1+x_2)(t+1)+1-(y_1+y_2) \\
   &=& a(t+1)+1-b \\
   &=& dt+d+r-d \\
   &=& dt+r=n,
\end{eqnarray*}
if we define the graph $E_6$ to be
 $$E_6: \ \ 
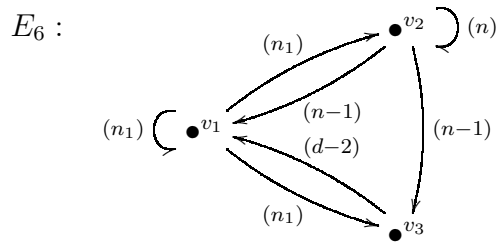
\xymatrix{
& & \bullet^{v_2} \ar@(ur,dr)[]^{(n)} \ar@/^/[dll]^{(n-1)} \ar@/^/[dd]^{(n-1)} \\
\bullet^{v_1} \ar@(ul,dl)[]_{(n_1)} \ar@/^/[urr]^{(n_1)}
\ar@/_/[drr]_{(n_1)} & & \\
 & & \bullet^{v_3}  \ar@/_/[ull]_{(d-2)}
 }$$
 then we get
 $E_5=E_6(v_1\hookrightarrow v_2),$ so that
 $$L(E_5)\cong L(E_6) \mbox{ by Theorem \ref{meu1_gen}.}$$

 \medskip

 {\bf Step 7.}   Define $E_7=E_6(v_3\hookrightarrow v_2).$  That is,
$$E_7: 
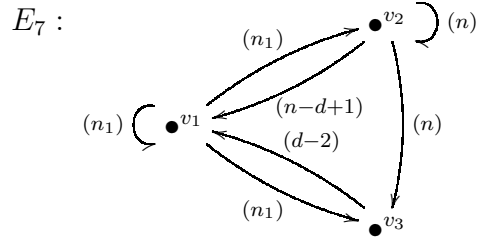
\xymatrix{
& & \bullet^{v_2} \ar@(ur,dr)[]^{(n)} \ar@/^/[dll]^{(n-d+1)} \ar@/^/[dd]^{(n)} \\
\bullet^{v_1} \ar@(ul,dl)[]_{(n_1)} \ar@/^/[urr]^{(n_1)}
\ar@/_/[drr]_{(n_1)} & & \\
 & & \bullet^{v_3}  \ar@/_/[ull]_{(d-2)}
 }$$
  Thus
  $$L(E_6)\cong L(E_7) \mbox{ by Theorem \ref{meu1_gen}.}$$

  \medskip

  {\bf Step 8.}   We let  $E_8$ denote the graph
  $$E_8: \ \ 
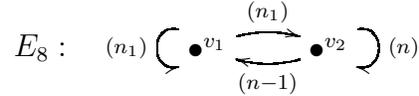
\xymatrix{
  \bullet^{v_1} \ar@(ul,dl)[]_{(n_1)} \ar@/^/[r]^{(n_1)} &
  \bullet^{v_2} \ar@(ur,dr)[]^{(n)} \ar@/^/[l]^{(n-1)} }$$
  We split the set of edges emitted by $v_2$ in two sets of $d-2$ and $2n-d+1$ edges respectively.  That is,
we use the partition
  $$
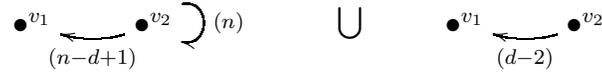
\xymatrix{
  \bullet^{v_1}  &
  \bullet^{v_2} \ar@(ur,dr)[]^{(n)} \ar@/^/[l]^{(n-d+1)} &  & \bigcup & \bullet^{v_1}  &
  \bullet^{v_2} \ar@/^/[l]^{(d-2)} }$$
Then it is not hard to see that $(E_8)_s(\emph{P})=E_7,$ so that
$$L(E_7)\cong L(E_8) \mbox{ by Theorem \ref{outsplitiso}.}$$

\medskip

{\bf Step 9.}   Now consider $E_9=E_8(v_1\hookrightarrow v_2).$  Pictorially,
$$E_9: \ \ \  
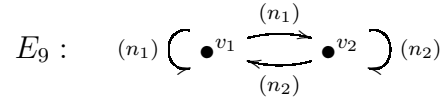
\xymatrix{
  \bullet^{v_1} \ar@(ul,dl)[]_{(n_1)} \ar@/^/[r]^{(n_1)} &
  \bullet^{v_2} \ar@(ur,dr)[]^{(n_2)} \ar@/^/[l]^{(n_2)} }$$
  But then
$$L(E_8)\cong L(E_9) \mbox{ by Theorem \ref{meu1_gen}.}$$

\medskip

{\bf Step 10.}   Finally, recall that the graph $R_n$ is given by

  $$R_n: \ \ \ \ 
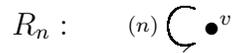
\xymatrix{
  \bullet^{v} \ar@(ul,dl)[]_{(n)}}$$
If we take a partition of the set of edges emitted by $v$ in two sets of $n_1$ and $n_2$ edges
respectively, then $(R_n)_s(\emph{P})=E_9,$
  so that
  $$L(E_9) \cong L(R_n) \mbox{ by Theorem \ref{outsplitiso}.}$$

\medskip

 Thus Steps 1 through 10 yield $L(R_n^d)\cong L(R_n)$, so that we have established the isomorphism
 $${\rm M}_d(L_n) \cong L(R_n^d) \cong L(R_n) \cong L_n,$$
 and we are done.
\end{proof}

\bigskip

We have seen above that Theorem \ref{meu1_gen} has a direct analog for C$^*$-algebras in Corollary
\ref{meu1_gen-C*}, and that Theorem \ref{outsplitiso} has a direct analog for C$^*$-algebras in
\cite[Theorem 3.2]{Flow}.  The proof of Theorem \ref{Th:Medium Fish} follows from Theorem
\ref{meu1_gen}, Theorem \ref{outsplitiso}, and from purely combinatorial arguments (arguments which
therefore hold in both the Leavitt path algebra and Cuntz-Krieger C$^*$-algebra settings).  So the
results of this section have provided a graph-theoretic approach to the following.

\begin{theor}\label{graphmatrixfishforCstaralgebras}
Let $n,d$ be positive integers having $\mbox{g.c.d.}(d,n-1)=1$. Then ${\mathcal
O}_n\cong\emph{M}_d({\mathcal O}_n).$
\end{theor}

As promised, we now show how Theorem \ref{Th:Medium Fish} yields the answer to another piece of The
Classification Question.  This same conclusion was drawn in
 \cite{AbAnhP};  for completeness, we present here some of the details of the proof provided there.

\begin{theor}\label{C*-fact2}\cite[Theorem 5.2]{AbAnhP}
Let $\mathcal{L}$ denote the set of purely infinite simple $K$-algebras
$$\{{\rm M}_d(L_n) | d,n \in \mathbb{N}\}.$$
Let $B,B'\in \mathcal{L}$. Then $B\cong B'$ if and only if there is an isomorphism $\phi:
K_0(B)\rightarrow K_0(B')$ for which $\phi([1_B])=[1_{B'}]$.
\end{theor}
\begin{proof}
It is well known that any unital isomorphism $f: B\rightarrow B'$ induces a group isomorphism
$K_0(f): K_0(B)\rightarrow K_0(B')$ sending $[1_B]$ to $[1_{B'}]$.

To see the converse, first notice that, for any $B\in \mathcal{L}$, $B={\rm M}_d(L_n)$ for suitable
$d, n\in \N$.
 As noted previously,
$(K_0({\rm M}_d(L_n)), [1_{{\rm M}_d(L_n)}])\cong (\Z/(n-1)\Z, \overline{d})$.  Hence, if $B'={\rm
M}_k(L_m)$ for suitable $k, m\in \N$, then the existence of an isomorphism $\phi: K_0(B)\rightarrow
K_0(B')$ forces that $n=m$.

Now, since every automorphism of $\Z/(n-1)\Z$ is given by multiplication by an element $1\leq l\leq
n-1$ such that ${\rm gcd}(l,n-1)=1$, the hypothesis $\phi([1_B])=[1_{B'}]$ yields that
$\overline{k}=\overline{dl}\in \Z/(n-1)\Z$, i.e., that $k\equiv dl$ (mod $n-1$). So Proposition
\ref{L:congruent mod n-1}
 gives that
$${\rm M}_k(L_n)\cong {\rm M}_{dl}(L_n)\cong {\rm M}_d({\rm M}_l(L_n)).$$
Since ${\rm g.c.d.}(l,n-1)=1$, we have ${\rm M}_l(L_n)\cong L_n$ by Theorem \ref{Th:Medium Fish}.
Hence, ${\rm M}_d({\rm M}_l(L_n))\cong {\rm M}_d(L_n)$, whence
$${\rm M}_k(L_n)\cong {\rm M}_{dl}(L_n)\cong {\rm M}_d({\rm M}_l(L_n))\cong {\rm M}_d(L_n),$$ as
desired.
\end{proof}

\bigskip

For reasons identical to those given prior to the statement of Theorem
\ref{graphmatrixfishforCstaralgebras}, Theorem \ref{C*-fact2} yields the following result for
 matrices over Cuntz-Krieger C$^*$-algebras.

 \begin{theor}
Let $\mathcal{L}$ denote the set of purely infinite simple {\rm C}$^*$-algebras
$$\{{\rm M}_d({\mathcal O}_n) | d,n \in \mathbb{N}\}.$$
Let $B,B'\in \mathcal{L}$. Then $B\cong B'$ if and only if there is an isomorphism $\phi:
K_0(B)\rightarrow K_0(B')$ for which $\phi([1_B])=[1_{B'}]$.
\end{theor}

\bigskip

We conclude with two remarks.  First, we observe that the sequence of isomorphisms between
$L(R_n^d)$ and $L(R_n)$ given in the proof of Theorem \ref{Th:Medium Fish} begins with a graph
having two vertices, eventually winds its way through graphs having three vertices, and finally
concludes with a graph having only one vertex. As it turns out, in any situation for which $t\geq
d-r-1$, we are in fact able to establish the isomorphism $L(R_n^d)$ and $L(R_n)$ without utilizing
graphs having three vertices. (However, we have been unable to achieve the desired isomorphism in
general without using graphs having three vertices.)   Second, we see no {\it a priori}  reason why
Theorems \ref{meu1_gen} and \ref{outsplitiso} should be expected to contain enough information to
provide us with proofs of The Classification Question in the situations described in both Sections
4 and 5.  That these two theorems suffice to cover both of these situations suggests that these two
theorems may indeed suffice to settle The Classification Question in further generality.

\end{document}